\documentclass{amsart}
\usepackage{amsmath, amssymb, latexsym}
\title{How to Construct a Flag Complex with a Given Face Vector}
\author{Andrew Frohmader}
\address{Department of Mathematics, 581 Malott Hall, Cornell University, Ithaca, NY 14853-4201}
\email{froh@math.cornell.edu}

\newtheorem{theorem}{Theorem}[section]
\newtheorem{proposition}[theorem]{Proposition}

\newtheorem{lemma}[theorem]{Lemma}
\newtheorem{definition}[theorem]{Definition}
\newtheorem{conjecture}[theorem]{Conjecture}
\newtheorem{example}[theorem]{Example}
\newtheorem{construction}[theorem]{Construction}

\def\proof{\smallskip\noindent {\it Proof: \ }}
\def\endproof{\hfill\ensuremath{\square}\medskip}

\begin{document}

\maketitle

\begin{abstract}
A method that often works for constructing a flag complex with a specified face vector is given.  This method can also be adapted to construct a vertex-decomposable (and hence Cohen-Macaulay) flag complex with a specified $h$-vector.
\end{abstract}

\section{Introduction}

There has been a lot of work done to characterize the face vectors of various classes of simplicial complexes.  Perhaps most famously, the Kruskal-Katona \cite{kruskal,katona} theorem characterizes the face vectors of all simplicial complexes.

In this paper, we are interested in flag complexes.  A simplicial complex is a flag complex if every minimal non-face is a two element set.  That is, given a set of vertices, either the entire set forms a face, or else some two of them do not form an edge.

Flag complexes are closely related to graphs.  The clique complex of a graph is a complex such that the vertices of the complex are the same as the vertices of the graph.  A set of vertices forms a face in the clique complex exactly if it forms a clique in the graph.  It is clear that if given a set of vertices in a graph, either they form a clique or else some two of them do not form an edge.  As such, the clique complex of a graph is a flag complex.  Conversely, every flag complex is the clique complex of its 1-skeleton, taken as a graph.

Some work has been done toward characterizing the face vectors of flag complexes.  A starting point is observing that every flag complex is, in particular, a simplicial complex.  As such, the bounds of the Kruskal-Katona theorem apply to flag complexes.

The converse of this does not work for flag complexes, however.  The Kruskal-Katona theorem also asserts that if a proposed integer vector satisfies certain bounds, then there is a simplicial complex that has the given vector as its face vector.  However, that simplicial complex may not be a flag complex, and there are many integer vectors that are the face vector of a simplicial complex, but not of a flag complex.

For example, (1, 3, 3) is the face vector of a simplicial complex, as one can take three vertices and assert that every pair of them forms an edge, while the three vertices do not form a triangle.  It is not the face vector of a flag complex, however, as in order to have three edges, every pair of vertices must form an edge, which forces the set of three vertices to form a triangle if the complex is to be flag.

There has been some work toward finding what sorts of integer vectors can be the face vectors of a simplicial complex, but not a flag complex.  For example, Lov\'{a}sz and Simonovits \cite{lovasz} found a lower bound on how many triangles a flag complex must have, if given its numbers of vertices and edges.  The example of the previous paragraph is a simple case of this.  For general simplicial complexes, there is nothing analogous to this, as however many vertices and edges a complex has, it can easily have no triangles whatsoever.

Another paper by this author \cite{mysecond} gave some bounds on the other side.  For example, the Kruskal-Katona theorem gives that a simplicial complex with nine triangles can have up to three faces of dimension three.  A flag complex with nine triangles can have at most two faces of dimension three.

These results both say that in order for an integer vector to be the face vector of a flag complex, it must satisfy certain bounds.  Equivalently, if it fails one of the bounds, then it is not the integer vector of a flag complex.

We do not have any non-trivial results that say, if an integer vector satisfies these certain conditions, then it must be the integer vector of a flag complex.  If given a particular integer vector that satisfies the known bounds, one can try to construct a flag complex with that particular vector as its face vector.  If one succeeds, then it is the face vector of a flag complex, but whether one succeeds can easily depend on how clever one is.  This is a wholly unsatisfactory situation, and one that this paper seeks to remedy.

Part of the problem is that flag complexes can behave strangely with regards to their face vectors.  Eckhoff \cite{eckhoff} showed that there are holes in what can be constructed.  For example, if a flag complex has exactly 15 edges, then it could have exactly 17 triangles or exactly 20 triangles, but it cannot have exactly 18 or 19 triangles.  If it has exactly 19 triangles, then it must have at least 17 edges.

In addition, bounds on consecutive face numbers are not enough.  \cite{mysecond}  For example, there is a flag complex with 70 triangles and 85 faces of dimension three.  There is also a flag complex with 85 faces of dimension three and 62 of dimension four.  However, these two complexes are rather different from each other, and a flag complex with 70 triangles can have at most 61 faces of dimension four, regardless of how many three-dimensional faces it has.

In Section~2, we start with the simplest case of considering only two face numbers at a time.  Construction~\ref{twofaceconst} often gives a flag complex with two specified face numbers.  Theorem~\ref{twofacetheorem} gives some bounds and shows that if the bounds are satisfied, then Construction~\ref{twofaceconst} works.  For example, there is a flag complex with exactly 7 faces of dimension five and 83 faces of dimension three.

Proposition~\ref{twofacebd4} gives an easily computed numerical bound that guarantees that the construction works.  This bound is pretty far from being sharp, so if two given face numbers satisfy the Kruskal-Katona bound but not our bound, the construction often still works.  Theorem~\ref{twofaceprop} gives an asymptotic bound on when the construction works.

The basic idea is that we give a particular way to attempt to construct a flag complex with the exact number of faces of the higher dimension that we want, and as few faces of the lower dimension as we can.  If we then want more faces of the lower dimension, it is easy to add a bunch of isolated cliques to add the lower dimensional faces without any new higher dimensional faces.

If we want to give a flag complex with a given face vector, however, we usually want to match the whole face vector, not just two particular face numbers.  A step in this direction is to not merely match those two face numbers, but to specify both the dimension of the complex and the two face numbers, and produce a flag complex both with the appropriate dimension and the two specified face numbers.  This is what we do in Section~3.

In this, a characterization of the face vectors of colored complexes is useful.  Frankl, F\"{u}redi, and Kalai \cite{balanced} gave some bounds analogous to those of the Kruskal-Katona theorem, and showed that the face vectors of a simplicial complex of a given chromatic number must satisfy their bounds.  Conversely, given an integer vector that does satisfy their bounds, they can provide a simplicial complex with that vector as its face vector and a chromatic number no greater than specified.

Another paper of this author \cite{myfirst} proved that the face vectors of flag complexes must satisfy these same bounds, if we use the size of the largest face of the flag complex in place of the chromatic number.  A simplicial complex of dimension $d-1$ has a set of $d$ vertices, all of which are adjacent to each other, so its chromatic number must be at least $d$.  The chromatic number could be far greater than $d$, however, so this sometimes means that the face vectors of flag complexes must satisfy much tighter bounds than those implied by their chromatic numbers.

This means that we have some numerical bounds on the face vectors of flag complexes of a given dimension.  Furthermore, the construction of Frankl, F\"{u}redi, and Kalai that attains their bound gives some clue on how to construct a flag complex that is close to attaining the bounds.  As with the case of the Kruskal-Katona theorem, however, the bounds only go one way for flag complexes.  That is, there are integer vectors that satisfy the bounds, but are not the face vector of any flag complex.

We use this in Construction~\ref{dimconst}, which is a slight modification of Construction~\ref{twofaceconst}.  This attempts to construct a flag complex of a given dimension with a given face number of one dimension, and as few faces as possible of a lower dimension.

As in the case of unrestricted dimension, Theorem~\ref{dimtheorem} gives bounds that characterize when the construction works.  Proposition~\ref{dimbd2} gives weaker numerical bounds that are easy to compute, and can guarantee that Construction~\ref{dimconst} works.  Theorem~\ref{dimprop} gives an asymptotic bound on when the construction works.

Next, in Section~4, we do some work to evaluate the bounds of the two constructions.  We compute some asymptotic bounds from those of the Kruskal-Katona and Frankl-F\"{u}redi-Kalai theorems. We show that the asymptotic bounds of Theorem~\ref{twofaceprop} are the best we could hope for with flag complexes, whether by Construction~\ref{twofaceconst} or any other.  The gap between Theorem~\ref{dimprop} and the bounds of the Frankl-F\"{u}redi-Kalai theorem is comparable to the gap between Theorem~\ref{twofaceprop} and the bounds of the Kruskal-Katona theorem, so Construction~\ref{dimconst} is likely to also be pretty good.

In Section~5, we give Construction~\ref{mainconst}, which is our main construction.  It attempts to construct a flag complex with an entire specified face vector, rather than only two particular face numbers.  It is very similar to applying Construction~\ref{dimconst} iteratively.  The basic idea of the construction is that we start by specifying the top dimensional faces, while using as few faces as we can of lower dimensions.  We then add any additional faces of the next dimension down that are necessary, while again using as few faces of smaller dimensions as possible.  We repeat this for faces of each smaller dimension until all of the faces required are in the complex.

This does not fully characterize the face vectors of flag complexes, however.  It leaves us with a situation where we have some bounds that we know that the face vector of a flag complex must satisfy, and a construction that often works for face vectors that satisfy those bounds, but sometimes fails.

In Section~6, we explain some tweaks that can sometimes make a construction work when Construction~\ref{mainconst} fails as written.  This does not give an explicit construction, but only some ways to modify Construction~\ref{mainconst}. Conjectures~\ref{bigconj} and~\ref{converseconj} would be great progress toward characterizing the face vectors of flag complexes.  They are nearly converses of each other.  We explain why there are compelling reasons to believe that the latter conjecture is true.  Conjecture~\ref{converseconj} is mainly of interest because it is nearly a converse to Conjecture~\ref{bigconj}.

Finally, in Section~7, Construction~\ref{cmconst} shows how to adapt Construction~\ref{mainconst} to produce a vertex-decomposable flag complex with a specified $h$-vector.  If Construction~\ref{mainconst} can produce a complex with a specified face vector, then Construction~\ref{cmconst} can modify it to produce a vertex-decomposable flag complex with exactly the same $h$-vector.

\section{Two face numbers}

Recall that a \textit{simplicial complex} $\Delta$ on a vertex set $W$ is a collection of subsets of $W$ such that (i) for every $v \in W$, $\{v\} \in \Delta$ and (ii) for every $B \in \Delta$, if $A \subset B$, then $A \in \Delta$.  The elements of $\Delta$ are called \textit{faces}.  A face on $i$ vertices is said to have \textit{dimension} $i-1$, while the dimension of a complex is the maximum dimension of a face of the complex. The \textit{$i$-th face number} of a simplicial complex $\Delta$, $f_{i-1}(\Delta)$, is the number of faces of $\Delta$ on $i$ vertices.  The \textit{face vector} of $\Delta$ lists the face numbers of $\Delta$.

The face vectors of simplicial complexes were characterized independently by Kruskal \cite{kruskal} and Katona \cite{katona}.  In order to state the Kruskal-Katona theorem, we need a lemma.

\begin{lemma} \label{kklemma}
Given positive integers $m$ and $k$, there is a unique way to pick integers $s \geq 0$ and $n_k, n_{k-1}, \dots, n_{k-s}$ such that
$$m = {n_k \choose k} + {n_{k-1} \choose k-1} + \dots + {n_{k-s} \choose k-s}$$
and $n_k > n_{k-1} > \dots > n_{k-s} \geq k-s > 0$.
\end{lemma}

The constants aren't that hard to compute, either.  We pick $n_k$ such that ${n_k \choose k} \leq m < {n_k + 1 \choose k}$, and then repeat with $k-1$ instead of $k$ and $m - {n_k \choose k}$ instead of $m$.  We stop when we would end up with 0 as the new value for $m$, and this is guaranteed to happen no later than when we use 1 for $k$.

\begin{theorem}[Kruskal-Katona] \label{kktheorem}
Let $\Delta$ be a simplicial complex and let $m = f_{k-1}(\Delta)$.  Let
$$m = {n_k \choose k} + {n_{k-1} \choose k-1} + \dots + {n_{k-s} \choose k-s}$$
as in Lemma~\ref{kklemma}.  Then
$$f_{k-2}(\Delta) \geq {n_k \choose k-1} + {n_{k-1} \choose k-2} + \dots + {n_{k-s} \choose k-s-1}.$$
Furthermore, if a positive integer vector with 1 as its first entry satisfies these inequalities for all $k \geq 1$, then it is the face vector of some simplicial complex.
\end{theorem}

The Kruskal-Katona theorem is sometimes stated in terms of upper bounds on the number of faces of one dimension higher, rather than lower bounds on the number of faces of one dimension lower.  This is equivalent to the formulation given above.

This also gives bounds on non-consecutive face numbers, as we can chain the bounds of the Kruskal-Katona theorem.  A known number of $f_{k-1}(\Delta)$ gives an upper bound on $f_k(\Delta)$, and then that value of $f_k(\Delta)$ gives an upper bound on $f_{k+1}(\Delta)$, and so forth.  More precisely, if
$$f_{k-1}(\Delta) = {n_k \choose k} + {n_{k-1} \choose k-1} + \dots + {n_{k-s} \choose k-s}$$
with the constants as in Lemma~\ref{kklemma}, then
$$f_{k-1-i}(\Delta) \geq {n_k \choose k-i} + {n_{k-1} \choose k-1-i} + \dots + {n_{k-s} \choose k-s-i}.$$
There are some technical details of what happens if $k-s-i < 0$, but the formula does work with the convention that ${n \choose k} = 0$ if $n > 0 > k$.

While the Kruskal-Katona theorem guarantees that there is a simplicial complex having the desired face vector, that complex doesn't have to be a flag complex.  If $s = 0$ or $s = 1$, then the standard ``rev-lex" complex does happen to be a flag complex, however.  If $s = 0$, we can take a clique complex of a clique on $n_k$ vertices.  If $s = 1$, we can take a clique complex of a graph that consists of a clique on $n_k$ vertices, plus one other vertex that is adjacent to $n_{k-1}$ of the previous vertices.  If $s \geq 2$, however, then the bounds of the Kruskal-Katona theorem are rarely attained by a flag complex; loosely, it only happens if $n_{k-2} - (k-2) \ll k$.

\begin{construction}  \label{twofaceconst}
\textup{Let $k$, $m$, $p$, and $q$ be positive integers with $k > p$.  We wish to create a flag complex $\Delta$ with $f_{k-1}(\Delta) = m$ and $f_{p-1}(\Delta) = q$.  Define $n_0$ to be the unique integer such that ${n_0 \choose k} \leq m < {n_0 + 1 \choose k}$.  Define $m_0 = m - {n_0 \choose k}$.}

\textup{Define $m_i$ and $n_i$ for $i \geq 1$ recursively as follows.  If $m_{i-1} = 0$, then let $z = i-1$ and stop.  If $m_{i-1} > 0$, then let $n_i$ be the unique integer such that ${n_i \choose k-1} \leq m_{i-1} < {n_i + 1 \choose k-1}$ and $m_i = m_{i-1} - {n_i \choose k-1}$.}

\textup{Construct a complex $\Delta$ as the clique complex of a graph that starts with a clique on $n_0$ vertices, and then adds $z$ additional vertices, $v_1, v_2, \dots, v_z$, such that $v_i$ is adjacent to $n_i$ vertices of the original clique.}

\textup{It is easy to compute that}
\begin{eqnarray*}
f_{k-1}(\Delta) & = & {n_0 \choose k} + {n_1 \choose k-1} + {n_2 \choose k-1} + \dots + {n_z \choose k-1} \textup{ and} \\ f_{p-1}(\Delta) & = & {n_0 \choose p} + {n_1 \choose p-1} + {n_2 \choose p-1} + \dots + {n_z \choose p-1}.
\end{eqnarray*}
\textup{It is also easy to see that}
$$m_{i-1} = {n_i \choose k-1} + {n_{i+1} \choose k-1} + \dots + {n_z \choose k-1},$$
\textup{by starting with $i = z$ and working backwards.  From this, it follows that $f_{k-1}(\Delta) = m$.}

\textup{If $f_{p-1}(\Delta) > q$, then the construction fails.  If $f_{p-1}(\Delta) \leq q$, then add $q - f_{p-1}(\Delta)$ vertices, each adjacent to $p-1$ vertices of the original clique on $n_0$ vertices.  This will make $f_{p-1}(\Delta) = q$ and leave $f_{k-1}(\Delta) = m$.  \endproof}
\end{construction}

There are a few comments necessary to make sense of the construction.  First, it will terminate in finitely many steps, as $m_i < m_{i-1}$, and both are integers.  Next, it is possible to add a new vertex adjacent to $n_i$ of the first $n_0$ vertices, as $n_i < n_0$, for otherwise, we would have
$$m \geq {n_0 \choose k} + {n_i \choose k-1} \geq {n_0 \choose k} + {n_0 \choose k-1} = {n_0 + 1 \choose k},$$
which contradicts the choice of $n_0$.  Also, it is possible to add a new vertex and make it adjacent to $p-1$ of the first $n_0$ vertices because $n_0 \geq k-1 > p-1$.

Furthermore, the way that each $n_i$ is chosen is well-defined. Because $k > p > 0$, we must have $k-1 \geq 1$.  Thus,
$$1 = {k-1 \choose k-1} < {k \choose k-1} < {k+1 \choose k-1} < \dots .$$
$m_{i-1}$ is a positive integer, so it must either be between some two consecutive terms of the chain of inequalities or equal to one of them.  Either way, the choice of $n_i$ is unique.  We can do the same for $n_0$, using $k$ rather than $k-1$.

The big question is whether, after adding the first $n_0 + z$ vertices, we will have $f_{p-1}(\Delta) \leq q$.  Obviously, in some cases we cannot, as the Kruskal-Katona theorem gives conditions under which there is no complex $\Delta$ with $f_{k-1}(\Delta) = m$ and $f_{p-1}(\Delta) = q$, whether flag or otherwise.  As such, we define some functions that basically parameterize when the construction works.

\begin{definition}
\textup{Let $k$ and $p$ be positive integers with $k > p$.  Define $g_k(m)$ for $m > 0$ to be the unique integer such that ${g_k(m) \choose k-1} \leq m < {g_k(m) + 1 \choose k-1}$.  Define $c_p^k(m)$ for $m \geq 0$ recursively by $c_p^k(0) = 0$ and $c_p^k(m) = {g_k(m) \choose p-1} + c_p^k \big(m - {g_k(m) \choose k-1} \big)$ for $m > 0$.}
\end{definition}

\begin{lemma} \label{twofaceconlemma}
Let $k$, $p$, and $q$ be nonnegative integers with $k > p > 0$, and let
$$q = {n_0 \choose k} + {n_1 \choose k-1} + {n_2 \choose k-1} + \dots + {n_z \choose k-1}$$
as in Construction~\ref{twofaceconst}.  Let
$$m = {n_2 \choose k-1} + \dots + {n_z \choose k-1}.$$
Then
$$c_p^k(m) = {n_2 \choose p-1} + \dots + {n_z \choose p-1}.$$
\end{lemma}

\proof  We use induction on $z$.  For the base case, if $z \leq 1$, then $m = 0$ and $c_p^k(0) = 0$.  For the inductive step, if the lemma holds for $z-1$, then we get
\begin{eqnarray*}
c_p^k(m) & = & {n_2 \choose p-1} + c_p^k \bigg(m - {n_2 \choose k-1} \bigg) \\ & = & {n_2 \choose p-1} + {n_3 \choose p-1} + \dots + {n_z \choose p-1}
\end{eqnarray*}
by using the definition of $c_p^k(m)$ on the first line and the inductive hypothesis to produce the second.  \endproof

\begin{lemma} \label{twofacelemma}
Let $k \geq 2$ and $q > 0$ be integers.  There are unique integers $a$, $b$, and $m$ such that $q = {a \choose k} + {b \choose k-1} + m$, $a > b$, and ${b \choose k-2} > m \geq 0$.
\end{lemma}

The proof of this is basically the same as the proof of Lemma~\ref{kklemma}.  We simply let $m = {n_{k-2} \choose k-2} + \dots + {n_{k-s} \choose k-s}$.  If $s \leq 1$, then $m = 0$. If $s = 0$, then $b = k-2$.

\begin{theorem} \label{twofacetheorem}
Let $k$, $p$, $v$, and $w$ be positive integers with $k > p$.  Let $v = {a \choose k} + {b \choose k-1} + m$ as in Lemma~\ref{twofacelemma}.  If $w \geq {a \choose p} + {b \choose p-1} + c_p^k(m)$, then Construction~\ref{twofaceconst} gives a flag complex $\Delta$ with $f_{k-1}(\Delta) = v$ and $f_{p-1}(\Delta) = w$.
\end{theorem}

\proof  From Lemma~\ref{twofaceconlemma}, it is clear that there is a complex with $f_{k-1}(\Delta) = v$ and $f_{p-1}(\Delta) \leq w$ at one point in the construction.  The construction then adds as many faces of dimension $p-1$ as needed without adding any larger faces.  \endproof

It may at first glance seem peculiar to pull out two terms rather than one, and make $m$ the sum of the remaining terms, as all but the first term is of the form ${n_i \choose k-1}$.  The reason for this is that the first two terms match those of the Kruskal-Katona theorem, and this far is attainable with a flag complex.  The question of how good of a construction this is is partially a question of how close to the bounds of the Kruskal-Katona theorem we can come.  For this, it is only the third and subsequent terms that differ from the Kruskal-Katona bounds.

Saying that we have a construction and a function defined to describe when the construction works isn't a terribly satisfactory solution to the problem, however.  It would be better to say how large the function tends to be, and how far from the bounds of the Kruskal-Katona theorem this is.

For this, we have two answers.  First is an easily computable upper bound on $c_p^k(m)$, as a function of $p$, $k$, and $m$. Second is an asymptotic answer of approximately how quickly $c_p^k(m)$ gets large for fixed $p$ and $k$ as $m$ becomes large.  As we will see in a later section, the asymptotic answer is about the best that we could possibly hope for.  First, we need a few lemmas.

\begin{lemma} \label{twofacebd1}
Let $s \geq k > p > 0$ be integers.  Then $s^{k-p}(s-k)^p \leq (s-p)^k$.
\end{lemma}

\proof  Applying the arithmetic mean-geometric mean inequality to a list of numbers consisting of $k-p$ copies of $s$ and $p$ copies of $s-k$ yields
\begin{eqnarray*}
{s(k-p) + p(s-k) \over k-p + p} & \geq & \sqrt[k-p+p]{s^{k-p}(s-k)^p} \\ s-p & \geq & \sqrt[k]{s^{k-p}(s-k)^p} \\ \hspace{41 mm} (s-p)^k & \geq & s^{k-p}(s-k)^p \hspace{41 mm} \square
\end{eqnarray*}

\begin{lemma} \label{twofacebd2}
Let $n \geq k > p > 0$ be integers.  Then
$${n \choose p}^k \leq {n \choose k}^p {k \choose p}^k.$$
\end{lemma}

\proof  We can write this out to get
$${(n!)^k \over (p!)^k((n-p)!)^k} \leq {(n!)^p \over (k!)^p((n-k)!)^p}{(k!)^k \over (p!)^k((k-p)!)^k}.$$
This simplifies to
$${(n!)^{k-p}((n-k)!)^p \over ((n-p)!)^k} \leq {(k!)^{k-p} \over ((k-p)!)^k}.$$
We show this by induction on $n$.  For the base case, if $n = k$, it is easy to see that equality holds.

For the inductive step, if the inequality holds and we replace $n$ by $n+1$, this does not change the right side, while it multiplies the left side by ${(n+1)^{k-p}(n+1-k)^p \over (n+1-p)^k}$.  If we let $s = n+1$ and apply Lemma~\ref{twofacebd1}, we get $(n+1)^{k-p}(n+1-k)^p \leq (n+1-p)^k$, from which ${(n+1)^{k-p}(n+1-k)^p \over (n+1-p)^k} \leq 1$.  Therefore, increasing $n$ by one leaves the decreases the left side of the inequality above while leaving the right side unchanged, so the inequality still holds.  \endproof

\begin{lemma} \label{twofacebd3}
Let $a > 0$, $b \geq 0$, and $0 < q < 1$ be real numbers.  Then $(a+b)^q - b^q \geq a(a+b)^{q-1}$.
\end{lemma}

\proof  Let $f(x) = x^q$ on $x > 0$.  The tangent line approximation to $f(x)$ at $x = a+b$ is $y = (a+b)^q + (x-a-b)(a+b)^{q-1}$.  Because $f''(x) = q(q-1)x^{q-2} < 0$, the curve is below the tangent line, and so $f(x) \leq (a+b)^q + (x-a-b)(a+b)^{q-1}$ for all $x \geq 0$.  Plug in $x = b$ to get $b^q \leq (a+b)^q - a(a+b)^{q-1}$, which rearranges to the statement of the lemma. \endproof

\begin{proposition} \label{twofacebd4}
Let $k > p > 1$ be integers.  Then
$$c_p^k(m) \leq {k-1 \choose p-1} \bigg({k+1 \over 2}\bigg)^{{k-p \over k-1}} m^{{p-1 \over k-1}}.$$
\end{proposition}

\proof  We use induction on $m$.  We let $n$ be the unique integer such that ${n \choose k-1} \leq m < {n+1 \choose k-1}$.  For the base case, if $m = 0$, the inequality is $c_p^k(0) \leq 0$, which is true because $c_p^k(0) = 0$.

For the inductive step, let $w = m - {n \choose k-1}$.  By the inductive hypothesis, the lemma holds for $w$, so
$$c_p^k(w) \leq {k-1 \choose p-1} \bigg({k+1 \over 2}\bigg)^{{k-p \over k-1}} w^{{p-1 \over k-1}}.$$
To show that it holds for $m$ is equivalent to showing
$$c_p^k\bigg(w+{n \choose k-1}\bigg) \leq {k-1 \choose p-1}^{k-1} \bigg(w+{n \choose k-1}\bigg)^{{p-1 \over k-1}}.$$

We start with the statement of Lemma~\ref{twofacebd2}, using $p-1$ and $k-1$ in place of $p$ and $k$, which is
$${n \choose p-1}^{k-1} \leq {n \choose k-1}^{p-1} {k-1 \choose p-1}^{k-1}.$$
Multiply both sides by $\big({n+1 \over n-k+2}\big)^{k-p}$ to get
$${n \choose p-1}^{k-1} \bigg({n+1 \over n-k+2}\bigg)^{k-p} \leq {n \choose k-1}^{p-1} {k-1 \choose p-1}^{k-1} \bigg({n+1 \over n-k+2}\bigg)^{k-p}.$$
We can multiply out some terms and rearrange to get
$${(n!)^{k-p} ((k-1)!)^{p-1} ((n-k+1)!)^{p-1} (n+1)^{k-p} \over ((p-1)!)^{k-1} ((n-p+1)!)^{k-1} (n-k+2)^{k-p}} \leq {k-1 \choose p-1}^{k-1} \bigg({n+1 \over n-k+2}\bigg)^{k-p}.$$
The left side rearranges as
$${(n!)^{k-1} ((k-1)!)^{k-1} ((n-k+1)!)^{k-1} ((n+1)!)^{k-p} \over ((p-1)!)^{k-1} ((n-p+1)!)^{k-1} (n!)^{k-1} ((k-1)!)^{k-p} ((n-k+2)!)^{k-p}},$$
which factors as
$${n \choose p-1}^{k-1} {n \choose k-1}^{1-k} {n+1 \choose k-1}^{k-p}.$$
Thus, we have
$${n \choose p-1}^{k-1} {n \choose k-1}^{1-k} {n+1 \choose k-1}^{k-p} \leq {k-1 \choose p-1}^{k-1} \bigg({n+1 \over n-k+2}\bigg)^{k-p},$$
or equivalently,
$${n \choose p-1}^{k-1} \leq {n \choose k-1}^{k-1} {k-1 \choose p-1}^{k-1} {n+1 \choose k-1}^{p-k} \bigg({n+1 \over n-k+2}\bigg)^{k-p}.$$
Since ${n+1 \over n-k+2}$ decreases as $n$ gets larger and $n \geq k$, ${n+1 \over n-k+2} \leq {k+1 \over 2}$.  This gives us
$${n \choose p-1}^{k-1} \leq {n \choose k-1}^{k-1} {k-1 \choose p-1}^{k-1} {n+1 \choose k-1}^{p-k} \bigg({k+1 \over 2}\bigg)^{k-p}.$$
Take the $(k-1)$-th root of both sides to get
$${n \choose p-1} \leq {n \choose k-1} {k-1 \choose p-1} {n+1 \choose k-1}^{{p-k \over k-1}} \bigg({k+1 \over 2}\bigg)^{{k-p \over k-1}}.$$
From the definition of $n$, we know that $w + {n \choose k-1} < {n+1 \choose k-1}$.  Since the exponent on ${n+1 \choose k-1}$ in the line above is negative,
$${n+1 \choose k-1}^{{p-k \over k-1}} \leq \bigg(w+{n \choose k-1}\bigg)^{{p-k \over k-1}}.$$
This gives us
$${n \choose p-1} \leq {n \choose k-1} {k-1 \choose p-1} \bigg({k+1 \over 2}\bigg)^{{k-p \over k-1}} \bigg(w+{n \choose k-1}\bigg)^{{p-k \over k-1}}.$$
If we apply Lemma~\ref{twofacebd3} with $a = {n \choose k-1}$, $b = w$, and $q = {p-1 \over k-1}$, then we get
$${n \choose k-1} \bigg(w+{n \choose k-1}\bigg)^{{p-1 \over k-1} - 1} \leq \bigg(w+{n \choose k-1}\bigg)^{{p-1 \over k-1}} - w^{{p-1 \over k-1}},$$
from which we have
$${n \choose p-1} \leq {k-1 \choose p-1} \bigg({k+1 \over 2}\bigg)^{{k-p \over k-1}} \bigg(\bigg(w+{n \choose k-1}\bigg)^{{p-1 \over k-1}} - w^{{p-1 \over k-1}}\bigg),$$
or equivalently
$${n \choose p-1} + w^{{p-1 \over k-1}}{k-1 \choose p-1} \bigg({k+1 \over 2}\bigg)^{{k-p \over k-1}} \leq {k-1 \choose p-1} \bigg({k+1 \over 2}\bigg)^{{k-p \over k-1}} \bigg(w+{n \choose k-1}\bigg)^{{p-1 \over k-1}}.$$
From the definition of $c_p^k(w)$, we know that
$$c_p^k\bigg(w+{n \choose k-1}\bigg) = {n \choose p-1} + c_p^k(w).$$
The inductive hypothesis gives us
$$c_p^k(w) \leq {k-1 \choose p-1} \bigg({k+1 \over 2}\bigg)^{{k-p \over k-1}} w^{{p-1 \over k-1}}.$$
Therefore,
\begin{eqnarray*}
c_p^k\bigg(w+{n \choose k-1}\bigg) & = & {n \choose p-1} + c_p^k(w) \\ & \leq & {n \choose p-1} + {k-1 \choose p-1} \bigg({k+1 \over 2}\bigg)^{{k-p \over k-1}} w^{{p-1 \over k-1}} \\ & \leq & {k-1 \choose p-1} \bigg({k+1 \over 2}\bigg)^{{k-p \over k-1}} \bigg(w+{n \choose k-1}\bigg)^{{p-1 \over k-1}}.
\end{eqnarray*}
This completes the inductive step, and hence the proof.  \endproof

We can use this lemma to prove an asymptotic bound on $c_p^k(m)$.

\begin{theorem}  \label{twofaceprop}
Let $k > p > 1$ be integers.  Then
$$\lim_{m \to \infty} {c_p^k(m) \over m^{{p-1 \over k-1}}} = {((k-1)!)^{{p-1 \over k-1}} \over (p-1)!}.$$
\end{theorem}

Note that this theorem both claims that the limit exists and also gives the value of the limit.

\proof  For a given integer $m > 0$, let $n(m)$ be the unique integer such that ${n(m) \choose k-1} \leq m < {n(m) + 1 \choose k-1}$.  We can compute ${(n(m) - k + 2)^{k-1} \over (k-1)!} \leq {n(m) \choose k-1} \leq m$, from which
$$m^{-{p-1 \over k-1}} \leq \bigg( {(n(m) - k + 2)^{k-1} \over (k-1)!} \bigg)^{-{p-1 \over k-1}} = (n(m) - k + 2)^{1-p} ((k-1)!)^{{p-1 \over k-1}}.$$

Let $w(m) = m - {n(m) \choose k-1}$.  We have $m = {n(m) \choose k-1} + w(m)$ with $w(m) < {n(m) \choose k-2}$, so by definition, $c_p^k(m) = {n(m) \choose p-1} + c_p^k(w(m))$.  By linearity,
$$\lim_{m \to \infty} {c_p^k(m) \over m^{{p-1 \over k-1}}} = \lim_{m \to \infty} {{n(m) \choose p-1} \over m^{{p-1 \over k-1}}} + \lim_{m \to \infty} {c_p^k(w(m)) \over m^{{p-1 \over k-1}}}$$
if both limits on the right hand side exist.

It is clear from the definition that $0 \leq c_p^k(w(m))$.  We can apply Proposition~\ref{twofacebd4} to $w(m)$ to get
\begin{eqnarray*}
& & {c_p^k(w(m)) \over m^{{p-1 \over k-1}}} \\ & \leq & m^{-{p-1 \over k-1}} {k-1 \choose p-1} \bigg({k+1 \over 2}\bigg)^{{k-p \over k-1}} w(m)^{{p-1 \over k-1}} \\ & \leq & (n(m) - k + 2)^{1-p} ((k-1)!)^{{p-1 \over k-1}} {k-1 \choose p-1} \bigg({k+1 \over 2}\bigg)^{{k-p \over k-1}} {n(m) \choose k-2}^{{p-1 \over k-1}} \\ & \leq & (n(m)-k+2)^{1-p} ((k-1)!)^{{p-1 \over k-1}} {k-1 \choose p-1} \bigg({k+1 \over 2}\bigg)^{{k-p \over k-1}} \bigg({(n(m))^{k-2} \over (k-2)!}\bigg)^{{p-1 \over k-1}} \\ & = & (n(m)-k+2)^{1-p} ((k-1)!)^{{p-1 \over k-1}} {k-1 \choose p-1} \bigg({k+1 \over 2}\bigg)^{{k-p \over k-1}} ((k-2)!)^{-{p-1 \over k-1}} n(m)^{{(k-2)(p-1) \over k-1}} \\ & = & (k-1)^{{p-1 \over k-1}} {k-1 \choose p-1} \bigg({k+1 \over 2}\bigg)^{{k-p \over k-1}} \bigg({n(m) \over n(m) - k + 2}\bigg)^{p-1} n(m)^{-{p-1 \over k-1}}.
\end{eqnarray*}
As $m \to \infty$, $n(m) \to \infty$ also, so ${n(m) \over n(m) - k + 2} \to 1$, from which $\big({n(m) \over n(m) - k + 2}\big)^{p-1} \to 1$.  In addition, $(k-1)^{{p-1 \over k-1}} {k-1 \choose p-1} \big({k+1 \over 2}\big)^{{k-p \over k-1}}$ is merely a messy constant that does not depend on $m$.  Finally, $n(m)^{-{p-1 \over k-1}} \to 0$.  The limit of the product is the product of the limits, so as $m \to \infty$, the final line goes to 0.  Since ${c_p^k(m) \over m^{{p-1 \over k-1}}}$ is non-negative and is bounded above by something that goes to zero as $m \to \infty$, $$\lim_{m \to \infty} {c_p^k(w(m)) \over m^{{p-1 \over k-1}}} = 0.$$

From ${n(m) \choose k-1} \leq m < {n(m) + 1 \choose k-1}$, we can compute
$${(n(m)-k+2)^{k-1} \over (k-1)!} \leq {n(m) \choose k-1} \leq m < {n(m) + 1 \choose k-1} \leq {(n(m)+1)^{k-1} \over (k-1)!}.$$
This allows us to compute
\begin{eqnarray*}
(n(m)-k+2)^{k-1} \leq & m(k-1)! & < (n(m)+1)^{k-1} \\ n(m)-k+2 \leq & (m(k-1)!)^{1 \over k-1} & < n(m)+1 \\ (m(k-1)!)^{1 \over k-1} - 1 < & n(m) & \leq (m(k-1)!)^{1 \over k-1} + k - 2.
\end{eqnarray*}
Since ${(n(m)-p+2)^{p-1} \over (p-1)!} \leq {n(m) \choose p-1} \leq {(n(m))^{p-1} \over (p-1)!}$, we can compute
\begin{eqnarray*}
{((m(k-1)!)^{1 \over k-1} - 1)^{p-1} \over (p-1)!} \leq & {n(m) \choose p-1} & \leq {((m(k-1)!)^{1 \over k-1} + k - 2)^{p-1} \over (p-1)!} \\ {((m(k-1)!)^{1 \over k-1} - 1)^{p-1} \over m^{{p-1 \over k-1}} (p-1)!} \leq & {{n(m) \choose p-1} \over m^{{p-1 \over k-1}}} & \leq {((m(k-1)!)^{1 \over k-1} + k - 2)^{p-1} \over m^{{p-1 \over k-1}} (p-1)!} \\ {(((k-1)!)^{1 \over k-1} - m^{{-1 \over k-1}})^{p-1} \over (p-1)!} \leq & {{n(m) \choose p-1} \over m^{{p-1 \over k-1}}} & \leq {(((k-1)!)^{1 \over k-1} + (k - 2)m^{{-1 \over k-1}})^{p-1} \over (p-1)!}.
\end{eqnarray*}
As $m \to \infty$, $m^{{-1 \over k-1}} \to 0$, so both ends of the last inequality go to ${((k-1)!)^{{p-1 \over k-1}} \over (p-1)!}$.  Therefore, so does the expression in the middle.  \endproof

As we will see in Section~4, this is the best asymptotic bound that we could possibly hope for.  We can use Stirling's approximation to get ${((k-1)!)^{{p-1 \over k-1}} \over (p-1)!} \approx \big({k-1 \over p-1}\big)^{p-{1 \over 2}}$.

\section{Restricting the dimension}

If we wish to construct a flag complex with two particular face numbers, then Construction~\ref{twofaceconst} works pretty well.  It will always have one very large face on as many vertices as possible, however.  If we have a particular dimension in mind for the flag complex, as is implicitly given if we have an entire face vector in mind, then the construction of the previous section probably won't have the suitable dimension.  Fortunately, there are some things known as to what can happen to a flag complex of a given dimension.

\begin{definition}
\textup{The Tur\'{a}n graph $T_{n,r}$ is the graph obtained by partitioning $n$ vertices into $r$ parts as evenly as possible, and making two vertices adjacent exactly if they are not in the same part.  We define ${n \choose k}_r$ to be the number of cliques on $k$ vertices in the Tur\'{a}n graph $T_{n,r}$.}
\end{definition}

Each part in the Tur\'{a}n graph has either $\big\lceil {n \over r} \big\rceil$ or $\big\lfloor {n \over r} \big\rfloor$ vertices.  If we define integers $p$ and $q$ by $n = pr + q$ with $0 < q \leq r$, then there are $q$ parts of $p+1$ vertices and $r-q$ parts of $p$ vertices.

There are a variety of formulas used for computing ${n \choose k}_r$ numerically, depending on how we wish to count cliques of the graph.  One useful identity is ${n \choose k}_r = {n-1 \choose k}_r + {n-p-1 \choose k-1}_{r-1}$, which is analogous the common combinatorial identity ${n \choose k} = {n-1 \choose k} + {n-1 \choose k-1}$, and follows from about the same proof of picking a vertex and counting the number of cliques that contain it and those that don't contain it.  One explicit formula is
$${n \choose k}_r = \sum_{i=0}^{q} {q \choose i}{r-i \choose k-i}p^{k-i}.$$
This follows from picking one ``last" vertex from each part with $p+1$ vertices, and letting $i$ be the number of ``last" vertices chosen.  It also works if we define $n = pr + q$ and require $0 \leq q < r$.

It is clear from the definition that if $pr \leq n \leq (p+1)r$, then ${pr \choose k}_r \leq {n \choose k}_r \leq {(p+1)r \choose k}_r$.  The formulas on the end are easy to compute, so we get the approximation ${r \choose k}p^r \leq {n \choose k}_r \leq {r \choose k}(p+1)^r$.

The face vectors of flag complexes are closely related to those of colored complexes.  A \textit{coloring} of a simplicial complex with color set $[n] = \{1, 2, \dots, n\}$ is an assignment of a color to each vertex of the complex such that no two vertices in the same face are the same color.  This is equivalent to no two vertices in the same edge being the same color, which is the same requirement as for a graph coloring of the 1-skeleton of the complex, taken as a graph.

A \textit{colored complex} is a simplicial complex with a coloring.  If a colored complex has $n$ colors, then it has dimension at most $n-1$, as an $n$ dimensional face would have $n+1$ vertices, no two of which would be the same color.  The face vectors of colored complexes were characterized by Frankl, F\"{u}redi, and Kalai. \cite{balanced}  To state their characterization, we need a lemma analogous to Lemma~\ref{kklemma}.

\begin{lemma} \label{colorcan}
Given positive integers m, k, and r with $r \geq k$, there are unique $s$, $n_k$, $n_{k-1}$, \dots , $n_{k-s}$ such that
$$m = {n_k\choose k}_r + {n_{k-1}\choose k-1}_{r-1} + \dots + {n_{k-s}\choose k-s}_{r-s},$$
$n_{k-i}-\big\lfloor{n_{k-i}\over r-i}\big\rfloor > n_{k-i-1}$ for all $0\leq i < s,$ and $n_{k-s}\geq k-s > 0$.
\end{lemma}

There is a larger gap then merely requiring $n_{k-i} > n_{k-i-1}$ for about the same reasons that the recursive identity is ${n \choose k}_r = {n-1 \choose k}_r + {n-p-1 \choose k-1}_{r-1}$ and not ${n \choose k}_r = {n-1 \choose k}_r + {n-1 \choose k-1}_{r-1}$.

\begin{theorem}[Frankl-F\"{u}redi-Kalai \cite{balanced}] \label{coloredkk}
Let $\Delta$ be an $r$-colored simplicial complex and let $m = f_{k-1}(\Delta)$.  Let
$$m = {n_k \choose k}_r + {n_{k-1} \choose k-1}_{r-1} + \dots + {n_{k-s} \choose k-s}_{r-s}$$
as in Lemma~\ref{colorcan}.  Then
$$f_{k-2}(\Delta) \geq {n_k \choose k-1}_r + {n_{k-1} \choose k-2}_{r-1} + \dots + {n_{k-s} \choose k-s-1}_{r-s}.$$
Furthermore, if a positive integer vector with 1 as its first entry satisfies these inequalities for all $k \geq 1$, then it is the face vector of some $r$-colored simplicial complex.
\end{theorem}

Flag complexes follow these same bounds.  \cite{myfirst}

\begin{theorem} \label{flagcolored}
Let $\Delta$ be a flag complex of dimension $d-1$.  Then there is a $d$-colored complex $\Gamma$ such that $f(\Gamma) = f(\Delta)$.
\end{theorem}

Because $\Gamma$ is $d$-colored, it has to follow the bounds of Theorem~\ref{coloredkk}.  This last theorem says that the flag complex $\Delta$ has to follow these same bounds, as though the complex were $d$-colorable.  If $\Delta$ has chromatic number $d$, then this is not surprising, but the chromatic number of $\Delta$ could be much larger than $d$.  This says that $\Delta$ must satisfy the bounds of Theorem~\ref{coloredkk} as though it had chromatic number $d$, even if its real chromatic number is much larger.

This theorem suggests that if we wish to construct a flag complex of a given dimension, we should look at Tur\'{a}n graphs.  As with the Kruskal-Katona theorem, we can attain the bounds of the theorem with a flag complex if $s \leq 1$, but usually not if $s$ is larger.  With this in mind, we offer a construction analogous to Construction~\ref{twofaceconst} that attains the bounds if $s \leq 1$ and tries to come close if $s \geq 2$.

\begin{construction}  \label{dimconst}
\textup{Let $k$, $m$, $p$, $q$, and $r$ be positive integers with $r \geq k > p$.  We wish to create a flag complex $\Delta$ of dimension $r-1$ with $f_{k-1}(\Delta) = m$ and $f_{p-1}(\Delta) = q$.  Define $n_0$ to be the unique integer such that ${n_0 \choose k}_r \leq m < {n_0 + 1 \choose k}_r$.  Define $m_0 = m - {n_0 \choose k}_r$.}

\textup{Define $m_i$ and $n_i$ for $i \geq 1$ recursively as follows.  If $m_{i-1} = 0$, then let $z = i-1$ and stop.  If $m_{i-1} > 0$, then let $n_i$ be the unique integer such that ${n_i \choose k-1}_{r-1} \leq m_{i-1} < {n_i + 1 \choose k-1}_{r-1}$ and $m_i = m_{i-1} - {n_i \choose k-1}_{r-1}$.}

\textup{Construct a complex $\Delta$ as the clique complex of a graph that starts with the Tur\'{a}n graph $T_{n_0,r}$, and then adds $z$ additional vertices, $v_1, v_2, \dots, v_z$, such that the link of $v_i$ is a Tur\'{a}n graph $T_{n_i,r-1}$ that is a subgraph of the original Tur\'{a}n graph $T_{n_0,r}$.}

\textup{It is easy to compute that}
\begin{eqnarray*}
f_{k-1}(\Delta) & = & {n_0 \choose k}_r + {n_1 \choose k-1}_{r-1} + {n_2 \choose k-1}_{r-1} + \dots + {n_z \choose k-1}_{r-1} \textup{ and} \\ f_{p-1}(\Delta) & = & {n_0 \choose p}_r + {n_1 \choose p-1}_{r-1} + {n_2 \choose p-1}_{r-1} + \dots + {n_z \choose p-1}_{r-1}.
\end{eqnarray*}
\textup{It is also easy to see that}
$$m_{i-1} = {n_i \choose k-1}_{r-1} + {n_{i+1} \choose k-1}_{r-1} + \dots + {n_z \choose k-1}_{r-1},$$
\textup{by starting with $i = z$ and working backwards.  From this, it follows that $f_{k-1}(\Delta) = m$.}

\textup{If $f_{p-1}(\Delta) > q$, then the construction fails.  If $f_{p-1}(\Delta) \leq q$, then add $q - f_{p-1}(\Delta)$ vertices, each adjacent to $p-1$ vertices of the original clique on $n_0$ vertices.  This will make $f_{p-1}(\Delta) = q$ and leave $f_{k-1}(\Delta) = m$.  \endproof}
\end{construction}

This construction makes sense for the same reasons as Construction~\ref{twofaceconst}.  As before, the big question is whether after adding the first $n_0 + z$ vertices, we will have $f_{p-1}(\Delta) \leq q$.  We define a function that parameterizes when the construction works.

\begin{definition}
\textup{Let $r$, $k$, and $p$ be positive integers with $r \geq k > p$.  Define $j_k^r(m)$ for $m > 0$ to be the unique integer such that ${j_k^r(m) \choose k-1}_{r-1} \leq m < {j_k^r(m) + 1 \choose k-1}_{r-1}$.  Define $d_{p,r}^k(m)$ for $m \geq 0$ recursively by $d_{p,r}^k(0) = 0$ and $d_{p,r}^k(m) = {j_k^r(m) \choose p-1} + d_{p,r}^k \big(m - {j_k^r(m) \choose k-1} \big)$ for $m > 0$.}
\end{definition}

\begin{lemma} \label{dimconlemma}
Let $r$, $k$, $p$, and $q$ be nonnegative integers with $r \geq k > p > 0$, and let
$$q = {n_0 \choose k}_r + {n_1 \choose k-1}_{r-1} + {n_2 \choose k-1}_{r-1} + \dots + {n_z \choose k-1}_{r-1}$$
as in Construction~\ref{dimconst}.  Let
$$m = {n_2 \choose k-1}_{r-1} + \dots + {n_z \choose k-1}_{r-1}.$$
Then
$$d_{p,r}^k(m) = {n_2 \choose p-1}_{r-1} + \dots + {n_z \choose p-1}_{r-1}.$$
\end{lemma}

\proof  We use induction on $z$.  For the base case, if $z \leq 1$, then $m = 0$ and $d_{p,r}^k(0) = 0$.  For the inductive step, if the lemma holds for $z-1$, then we get
\begin{eqnarray*}
d_{p,r}^k(m) & = & {n_2 \choose p-1}_{r-1} + d_{p,r}^k \bigg(m - {n_2 \choose k-1}_{r-1} \bigg) \\ & = & {n_2 \choose p-1}_{r-1} + {n_3 \choose p-1}_{r-1} + \dots + {n_z \choose p-1}_{r-1}
\end{eqnarray*}
by using the definition of $d_{p,r}^k(m)$ on the first line and the inductive hypothesis to produce the second.  \endproof

\begin{lemma} \label{dimlemma}
Let $r \geq k \geq 2$ and $q > 0$ be integers.  There are unique integers $a$, $b$, and $m$ such that $q = {a \choose k}_r + {b \choose k-1}_{r-1} + m$, $a - \big\lfloor {a \over r} \big\rfloor > b$, and ${b - \lfloor {b \over r-1} \rfloor \choose k-2}_{r-2} > m \geq 0$.
\end{lemma}

This lemma follows from Lemma~\ref{colorcan} in the same way that Lemma~\ref{twofacelemma} follows from Lemma~\ref{kklemma}.

\begin{theorem} \label{dimtheorem}
Let $r$, $k$, $p$, $v$, and $w$ be positive integers with $r \geq k > p$.  Let $v = {a \choose k}_r + {b \choose k-1}_{r-1} + m$ as in Lemma~\ref{dimlemma}.  If $w \geq {a \choose p}_r + {b \choose p-1}_{r-1} + d_{p,r}^k(m)$, then Construction~\ref{dimconst} gives a flag complex $\Delta$ of dimension at most $r-1$ with $f_{k-1}(\Delta) = v$ and $f_{p-1}(\Delta) = w$.
\end{theorem}

\proof  From Lemma~\ref{dimconlemma}, it is clear that there is a complex with $f_{k-1}(\Delta) = v$ and $f_{p-1}(\Delta) \leq w$ at one point in the construction.  The construction then adds as many faces of dimension $p-1$ as needed without adding any larger faces.  \endproof

As in the previous section, we wish to evaluate how quickly $d_{p,r}^k(m)$ grows as $m$ does.

\begin{lemma} \label{dimbd1}
Let $r \geq k > p > 1$ be integers.  Then
$${k-1 \choose p-1} \bigg({k+1 \over 2}\bigg)^{{k-p \over k-1}} \leq {r-1 \choose p-1} {r-1 \choose k-1}^{-{p-1 \over k-1}} 2^{k-1}.$$
\end{lemma}

\proof  The right-hand side can be written as
$${(r-1)! ((r-k)!)^{{p-1 \over k-1}} ((k-1)!)^{{p-1 \over k-1}} \over (r-p)! (p-1)! ((r-1)!)^{{p-1 \over k-1}}} 2^{k-1} = {((r-1)!)^{{k-p \over k-1}} ((r-k)!)^{{p-1 \over k-1}} ((k-1)!)^{{p-1 \over k-1}} \over (r-p)! (p-1)!} 2^{k-1}.$$
If we increase the value of $r$ by 1, this multiplies the right-hand side by
$${r^{{k-p \over k-1}} (r-k+1)^{{p-1 \over k-1}} \over r-p+1}.$$
By Lemma~\ref{twofacebd1}, if we use $r$ for $s$, $p-1$ for $p$, and $k-1$ for $k$, we get $r^{k-p} (r-k+1)^{p-1} \leq (r-p+1)^{k-1}$.  Divide both sides by $(r-p+1)^{k-1}$ and then take the $k$-th root of both sides to get
$${r^{{k-p \over k-1}} (r-k+1)^{{p-1 \over k-1}} \over r-p+1} \leq 1.$$
Therefore, increasing $r$ decreases the right-hand side of the inequality of the lemma.  As such, proving it for one value of $r$ proves the lemma for all smaller values of $r$.  It thus suffices to prove the lemma for the limit as $r \to \infty$ of the right-hand side of the lemma.

The right-hand side of the inequality of the lemma can be rewritten as
$${((r-1)(r-2) \dots (r-k+1))^{{k-p \over k-1}} (r-k)! ((k-1)!)^{{p-1 \over k-1}} \over (r-p)(r-p-1) \dots (r-k+1) (r-k)! (p-1)!} 2^{k-1}.$$
The $(r-k)!$ terms cancel, and we can bound the rest to get
$$\bigg({r-k+1 \over r-p}\bigg)^{k-p} {((k-1)!)^{{p-1 \over k-1}} \over (p-1)!} \leq {r-1 \choose p-1} {r-1 \choose k-1}^{-{p-1 \over k-1}} \leq \bigg({r-1 \over r-k+1}\bigg)^{k-p} {((k-1)!)^{{p-1 \over k-1}} \over (p-1)!}.$$
As $r \to \infty$, both ${r-k+1 \over r-p} \to 1$ and ${r-1 \over r-k+1} \to 1$, so the left and right ends of this chain of inequalities both go to ${((k-1)!)^{{p-1 \over k-1}} \over (p-1)!}$.  Therefore, so does the term in the middle, and we can multiply by $2^{k-1}$ to get
$$\lim_{r \to \infty} {r-1 \choose p-1} {r-1 \choose k-1}^{-{p-1 \over k-1}} 2^{k-1} = {((k-1)!)^{{p-1 \over k-1}} \over (p-1)!} 2^{k-1}.$$

Now we have eliminated $r$, and it suffices to show that
\begin{eqnarray*}
{k-1 \choose p-1} \bigg({k+1 \over 2}\bigg)^{{k-p \over k-1}} & \leq & {((k-1)!)^{{p-1 \over k-1}} \over (p-1)!} 2^{k-1} \\ {(k-1)! \over (k-p)! (p-1)!} \bigg({k+1 \over 2}\bigg)^{{k-p \over k-1}} & \leq & {((k-1)!)^{{p-1 \over k-1}} \over (p-1)!} 2^{k-1} \\ ((k-1)!)^{{k-p \over k-1}} \bigg({k+1 \over 2}\bigg)^{{k-p \over k-1}} & \leq & (k-p)! 2^{k-1} \\ ((k-1)!)^{k-p} \bigg({k+1 \over 2}\bigg)^{k-p} & \leq & ((k-p)!)^{k-1} 2^{(k-1)^2}.
\end{eqnarray*}
By Stirling's approximation, $(k-1)! < \sqrt{2\pi(k-1)} \big({k-1 \over e}\big)^{k-1} e^{1 \over 12(k-1)}$ and $(k-p)! > \sqrt{2\pi(k-p)} \big({k-p \over e}\big)^{k-p}$.  Thus, it suffices to show
$$\bigg(\sqrt{2\pi(k-1)} \bigg({k-1 \over e}\bigg)^{k-1} e^{1 \over 12(k-1)} \bigg)^{k-p} \bigg({k+1 \over 2}\bigg)^{k-p} \leq \bigg(\sqrt{2\pi(k-p)} \bigg({k-p \over e}\bigg)^{k-p} \bigg)^{k-1} 2^{(k-1)^2},$$
which simplifies to
$$\bigg({k-1 \over k-p}\bigg)^{(k-1)(k-p)} \sqrt{{(k-1)^{k-p} \over (k-p)^{k-1}}} e^{k-p \over 12(k-1)} \bigg({k+1 \over 2}\bigg)^{k-p} \leq (\sqrt{2\pi})^{p-1} 2^{(k-1)^2}.$$

We wish to fix $k$ and find an upper bound on $\big({k-1 \over k-p}\big)^{k-p}$ subject to $1 \leq p \leq k-1$.  If we drop the requirement that $p$ must be an integer, we can take the derivative with respect to $p$ and get
$${\partial \over \partial p} \bigg({k-1 \over k-p}\bigg)^{k-p} = \bigg({k-1 \over k-p}\bigg)^{k-p} \bigg(1 - \ln\bigg({k-1 \over k-p}\bigg)\bigg).$$
This is zero precisely when $1 - \ln\big({k-1 \over k-p}\big) = 0$, which happens only when ${k-1 \over k-p} = e$.  In this case, $\big({k-1 \over k-p}\big)^{k-p} = \big(e^{1 \over e}\big)^{k-1}$.  We also check the endpoints and get that if $p=1$, then $\big({k-1 \over k-p}\big)^{k-p} = 1$ and if $p=k-1$, then $\big({k-1 \over k-p}\big)^{k-p} = k-1$.  Since $k \geq 3$, it is easy to show that $\big(e^{1 \over e}\big)^{k-1} > k-1$.  Thus, $\big({k-1 \over k-p}\big)^{k-p} \leq \big(e^{1 \over e}\big)^{k-1}$, and so $\big({k-1 \over k-p}\big)^{(k-1)(k-p)} \leq \big(e^{1 \over e}\big)^{(k-1)^2}$.

Next, we wish to show that
$$\sqrt{{(k-1)^{k-p} \over (k-p)^{k-1}}} e^{k-p \over 12(k-1)} \bigg({k+1 \over 2}\bigg)^{k-p} \leq (\sqrt{2\pi})^{p-1} \bigg(2 e^{-{1 \over e}}\bigg)^{(k-1)^2}.$$
Since $p$ and $k$ are integers, this is easily verified by brute force for $k \leq 12$.  If $k \geq 13$, then it is easy to show that $(k-1)^{{1 \over k-1}} \leq \big(2 e^{-{1 \over e}}\big)^{{2 \over 3}}$, or equivalently, $(k-1)^{{3 \over 2}(k-1)} \leq \big(2 e^{-{1 \over e}}\big)^{(k-1)^2}$.  It is also easy to see that ${k+1 \over 2} \leq k-1$.  We can also compute $ e^{k-p \over 12(k-1)} \leq e^{1 \over 12} < \sqrt{2\pi} \leq (\sqrt{2\pi})^{p-1}$.  From these, we can compute
\begin{eqnarray*}
\sqrt{{(k-1)^{k-p} \over (k-p)^{k-1}}} e^{k-p \over 12(k-1)} \bigg({k+1 \over 2}\bigg)^{k-p} & \leq & (k-1)^{{k-1 \over 2}} (\sqrt{2\pi})^{p-1} (k-1)^{k-1} \\ & \leq & (\sqrt{2\pi})^{p-1} \bigg(2 e^{-{1 \over e}}\bigg)^{(k-1)^2}.
\end{eqnarray*}
We multiply this last line by the inequality $\big({k-1 \over k-p}\big)^{(k-1)(k-p)} \leq \big(e^{1 \over e}\big)^{(k-1)^2}$ to get
$$\bigg({k-1 \over k-p}\bigg)^{(k-1)(k-p)} \sqrt{{(k-1)^{k-p} \over (k-p)^{k-1}}} e^{k-p \over 12(k-1)} \bigg({k+1 \over 2}\bigg)^{k-p} \leq (\sqrt{2\pi})^{k-p} 2^{(k-1)^2},$$
which we showed earlier was sufficient to prove the lemma.  \endproof

\begin{proposition} \label{dimbd2}
Let $r \geq k > p > 1$ be integers.  Then
$$d_{p,r}^k(m) \leq {r-1 \choose p-1}{r-1 \choose k-1}^{-{p-1 \over k-1}} 2^{k-1} m^{{p-1 \over k-1}}.$$
\end{proposition}

We use induction on $m$.  For the base case, if $m < {r-1 \choose k-1}$, then each time we add a new vertex in Construction~\ref{dimconst}, all of its vertices are adjacent to each other.  Thus, we add exactly as many faces of all dimensions as we would have done in Construction~\ref{twofaceconst}.  As such, $d_{p,r}^k(m) = c_p^k(m)$.  We can apply Proposition~\ref{twofacebd4} and Lemma~\ref{dimbd1} in order to get
$$d_{p,r}^k(m) = c_p^k(m) \leq {k-1 \choose p-1} \bigg({k+1 \over 2}\bigg)^{{k-p \over k-1}} m^{{p-1 \over k-1}} \leq {r-1 \choose p-1} {r-1 \choose k-1}^{-{p-1 \over k-1}} 2^{k-1} m^{{p-1 \over k-1}}.$$

Otherwise, $m \geq {r-1 \choose k-1}$.  Let $n$ be the unique integer such that ${n \choose k-1}_{r-1} \leq m < {n+1 \choose k-1}_{r-1}$, let $w = m - {n \choose k-1}_{r-1}$, and let $s = \big\lfloor {n \over r-1} \big\rfloor$.  Note that $m \geq {r-1 \choose k-1}$ means $n \geq r-1$, and so $s \geq 1$.  We can compute
\begin{eqnarray*}
& & d_{p,r}^k(m) \\ & = & d_{p,r}^k\bigg({n \choose k-1}_{r-1} + w\bigg) \\ & = & {n \choose p-1}_{r-1} + d_{p,r}^k(w) \\ & \leq & {(p-1)(s+1) \choose p-1}_{r-1}+ d_{p,r}^k(w) \\ & = & {r-1 \choose p-1}(s+1)^{p-1}+ d_{p,r}^k(w) \\ & \leq & {r-1 \choose p-1}(s+1)^{p-1} \bigg({2s \over s+1}\bigg)^{k-1} + d_{p,r}^k(w) \\ & = & {r-1 \choose p-1}(s+1)^{p-k} 2^{k-1} s^{k-1} + d_{p,r}^k(w) \\ & = & {r-1 \choose p-1} {r-1 \choose k-1}^{-{p-1 \over k-1}} 2^{k-1} {r-1 \choose k-1}^{{p-k \over k-1}} (s+1)^{p-k} {r-1 \choose k-1}s^{k-1} + d_{p,r}^k(w) \\ & = & {r-1 \choose p-1} {r-1 \choose k-1}^{-{p-1 \over k-1}} 2^{k-1} {(s+1)(r-1) \choose k-1}_{r-1}^{{p-k \over k-1}} {s(r-1) \choose k-1}_{r-1} + d_{p,r}^k(w) \\ & \leq & {r-1 \choose p-1} {r-1 \choose k-1}^{-{p-1 \over k-1}} 2^{k-1} {n+1 \choose k-1}_{r-1}^{{p-k \over k-1}} {n \choose k-1}_{r-1} + d_{p,r}^k(w) \\ & \leq & {r-1 \choose p-1} {r-1 \choose k-1}^{-{p-1 \over k-1}} 2^{k-1} m^{{p-k \over k-1}} {n \choose k-1}_{r-1} + d_{p,r}^k(w) \\ & \leq & {r-1 \choose p-1} {r-1 \choose k-1}^{-{p-1 \over k-1}} 2^{k-1} \bigg(m^{{p-1 \over k-1}-1}{n \choose k-1}_{r-1} + w^{{p-1 \over k-1}} \bigg) \\ & \leq & {r-1 \choose p-1}{r-1 \choose k-1}^{-{p-1 \over k-1}} 2^{k-1} m^{{p-1 \over k-1}}.
\end{eqnarray*}
The tenth and eleventh lines use that $k > p$, so that ${p-k \over k-1} < 0$, and making $c > 0$ smaller makes $c^{p-k \over k-1}$ larger.  The twelfth line uses the inductive hypothesis, and the final line uses Lemma~\ref{twofacebd3} with $a = {n \choose k-1}_{r-1}$, $b = w$, and $q = {p-1 \over k-1}$.  \endproof

\begin{theorem} \label{dimprop}
Let $r \geq k > p > 1$ be integers.
$$\lim_{m \to \infty} {d_{p,r}^k(m) \over m^{{p-1 \over k-1}}} = {r-1 \choose p-1}{r-1 \choose k-1}^{-{p-1 \over k-1}}.$$
\end{theorem}

As with Theorem~\ref{twofaceprop}, this both asserts that the limit exists and gives the value.

\proof  Let $n(m)$ be the unique integer such that ${n(m) \choose k-1}_{r-1} \leq m < {n(m)+1 \choose k-1}_{r-1}$, let $w(m) = m - {n(m) \choose k-1}_{r-1}$, and let $s(m) = \big\lfloor {n(m) \over r-1} \big\rfloor$.  Since ${n(m) \choose k-1}_{r-1} + {n(m) - s(m) \choose k-2}_{r-2} = {n(m) + 1 \choose k-1}_{r-1}$, we get $w(m) < {n(m) - s(m) \choose k-2}_{r-2} \leq {r-2 \choose k-2}(s+1)^{k-2}$.  For simplicity, let $\alpha = {r-1 \choose p-1}{r-1 \choose k-1}^{-{p-1 \over k-1}}$.

By definition, $d_{p,r}^k(m) = {n(m) \choose p-1}_{r-1} + d_{p,r}^k(w(m))$.  We can use this and Proposition~\ref{dimbd2} to get
\begin{eqnarray*}
{{n(m) \choose p-1}_{r-1} \over m^{{p-1 \over k-1}}} \leq & {d_{p,r}^k(m) \over m^{{p-1 \over k-1}}} & \leq {{n(m) \choose p-1}_{r-1} + d_{p,r}^k(w(m)) \over m^{{p-1 \over k-1}}} \\ {{n(m) \choose p-1}_{r-1} \over {n(m)+1 \choose k-1}_{r-1}^{{p-1 \over k-1}}} \leq & {d_{p,r}^k(m) \over m^{{p-1 \over k-1}}} & \leq {{n(m) \choose p-1}_{r-1} + \alpha 2^{k-1} w(m)^{{p-1 \over k-1}} \over {n(m) \choose k-1}_{r-1}^{{p-1 \over k-1}}} \\ {{r-1 \choose p-1}s(m)^{p-1} \over \big({r-1 \choose k-1}(s(m)+1)^{k-1}\big)^{{p-1 \over k-1}}} \leq & {d_{p,r}^k(m) \over m^{{p-1 \over k-1}}} & \leq {{r-1 \choose p-1}(s(m)+1)^{p-1} + \alpha 2^{k-1} \big({r-2 \choose k-2}(s(m)+1)^{k-2}\big)^{{p-1 \over k-1}} \over \big({r-1 \choose k-1}s(m)^{k-1}\big)^{{p-1 \over k-1}}} \\ \alpha \bigg({s(m) \over s(m)+1}\bigg)^{p-1} \leq & {d_{p,r}^k(m) \over m^{{p-1 \over k-1}}} & \leq \alpha \bigg({s(m)+1 \over s(m)}\bigg)^{p-1} + \alpha 2^{k-1} \bigg({{r-2 \choose k-2} \over {r-1 \choose k-1}} \bigg({s(m)+1 \over s(m)}\bigg)^{k-2} {1 \over s(m)} \bigg)^{{p-1 \over k-1}}.
\end{eqnarray*}
As $m \to \infty$, $s(m) \to \infty$.  Thus, ${s(m) \over s(m)+1} \to 1$ and ${s(m)+1 \over s(m)} \to 1$.  Therefore, the left-hand side and the first term on the right-hand side both go to $\alpha$ as $m \to \infty$.  The second term on the right-hand side goes to zero, as ${1 \over s(m)} \to 0$, while everything else in the second term goes to some finite constant.  Since ${d_{p,r}^k(m) \over m^{{p-1 \over k-1}}}$ is bounded between two formulas that both go to $\alpha$ as $m \to \infty$, so does ${d_{p,r}^k(m) \over m^{{p-1 \over k-1}}}$.  \endproof

Recall that
$$\lim_{r \to \infty} {r-1 \choose p-1}{r-1 \choose k-1}^{-{p-1 \over k-1}} = {((k-1)!)^{{p-1 \over k-1}} \over (p-1)!},$$
as shown in the proof of Lemma~\ref{dimbd1}.  This is perhaps to be expected, as one could think of the case where there are no color restrictions as being the limit as the number of available colors goes to infinity.

\section{Some analysis and improvements}

In this section, we compare the limits of Theorems~\ref{twofaceprop} and~\ref{dimprop} to those resulting from the bounds of the Kruskal-Katona and Frankl-F\"{u}redi-Kalai theorems.  We then argue that the limits of Theorems~\ref{twofaceprop} and~\ref{dimprop} are pretty good.  Finally, we give an improvement to Constructions~\ref{twofaceconst} and ~\ref{dimconst} that we will make use of in the next section, but does not change the limits of Theorems~\ref{twofaceprop} and~\ref{dimprop}.

First, we want to get a limit for the Kruskal-Katona theorem analogous to that of Theorem~\ref{twofaceprop}.

\begin{definition}
Let $m > 0$ and $k > p \geq 1$ be integers.  Let
$$m = {n_k \choose k} + {n_{k-1} \choose k-1} + \dots + {n_{k-s} \choose k-s}$$
as in Lemma~\ref{kklemma}.  Define
$$g_p^k(m) = {n_k \choose p} + {n_{k-1} \choose p-1} + \dots + {n_{k-s} \choose p-s}.$$
\end{definition}

With this definition, Theorem~\ref{kktheorem} essentially states that if $\Delta$ is a simplicial complex with $f_{k-1}(\Delta) = m$, then $f_{p-1}(\Delta) \geq g_p^k(m)$.

\begin{lemma} \label{kklimit}
Let $k > p \geq 1$ be integers and let $f(x)$ be a function whose domain is the natural numbers such that for positive integers $m$ and $n$, if ${n \choose k} \leq m \leq {n+1 \choose k}$, then ${n \choose p} \leq f(m) \leq {n+1 \choose p}$.  Then
$$\lim_{m \to \infty} {f(m) \over m^{p \over k}} = {(k!)^{{p \over k}} \over p!}.$$
\end{lemma}

The use of this lemma is that $g_p^k(m)$ is a function that satisfies the requirements.

\proof  For each integer $m$, let $n(m)$ be the unique integer such that ${n(m) \choose k} \leq m < {n(m) + 1 \choose k}$.  We can compute
\begin{eqnarray*}
{{n(m) \choose p} \over {n(m)+1 \choose k}^{p \over k}} \leq & {f(m) \over m^{p \over k}} & \leq {{n(m)+1 \choose p} \over {n(m) \choose k}^{p \over k}} \\ {{(n(m)-p)^p \over p!} \over \big({(n(m)+1)^k \over k!}\big)^{p \over k}} \leq & {f(m) \over m^{p \over k}} & \leq {{(n(m)+1)^p \over p!} \over \big({(n(m)-k)^k \over k!}\big)^{p \over k}} \\ \bigg( {n(m) - p \over n(m) + 1} \bigg)^p {(k!)^{{p \over k}} \over p!} \leq & {f(m) \over m^{p \over k}} & \leq \bigg( {n(m) + 1 \over n(m) - k} \bigg)^p {(k!)^{{p \over k}} \over p!}.
\end{eqnarray*}
As $m \to \infty$, $n(m) \to \infty$ also.  Thus, ${n(m) - p \over n(m) + 1} \to 1$ and ${n(m) + 1 \over n(m) - k} \to 1$.  As such, both ends of the chain of inequalities go to ${(k!)^{{p \over k}} \over p!}$ as $m \to \infty$.  Since ${f(m) \over m^{p \over k}}$ is bounded between two things that both go to ${(k!)^{{p \over k}} \over p!}$, it does as well.  \endproof

Lemma~\ref{kklimit} immediately gives us that
$$\lim_{m \to \infty} {g_p^k(m) \over m^{p \over k}} = {(k!)^{{p \over k}} \over p!}.$$
Using Proposition~\ref{twofacebd4}, it is not difficult to make a proof analogous to that of Theorem~\ref{twofaceprop} to show that if we have $q = {a(q) \choose k} + {b(q) \choose k-1} + m(q)$ as in Lemma~\ref{twofacelemma}, then
$$\lim_{q \to \infty} {{a(q) \choose p} + {b(q) \choose p-1} + c_p^k(m(q)) \over \big({a(q) \choose k} + {b(q) \choose k-1} + m(q)\big)^{{p-1 \over k-1}}} = {(k!)^{{p \over k}} \over p!}.$$
This shouldn't be all that surprising, considering that both the standard construction that attains the bound of the Kruskal-Katona theorem and Construction~\ref{twofaceconst} start with a simplex on $a(q)$ vertices.  Likewise, if we subtract off the first term from both constructions, what we get the same asymptotic limits for both constructions.

It is only in the third term that the constructions differ.  Here, if we chop off the first two terms of the Kruskal-Katona bound, we are left with $m = {n_{k-2} \choose k-2} + \dots + {n_{k-s} \choose k-s}$ and we can compute
$$\lim_{m \to \infty} {g_{p-2}^{k-2}(m) \over m^{{p-2 \over k-2}}} = {((k-2)!)^{{p-2 \over k-2}} \over (p-2)!}.$$
This differs from the bound of Theorem~\ref{twofaceprop}, which asserted that
$$\lim_{m \to \infty} {c_p^k(m) \over m^{{p-1 \over k-1}}} = {((k-1)!)^{{p-1 \over k-1}} \over (p-1)!}.$$
Because the power of $m$ in the denominator is different, we immediately get
$$\lim_{m \to \infty} {g_{p-2}^{k-2}(m) \over m^{{p-1 \over k-1}}} = 0.$$

At first glance, one might think that this means that there is a lot of room to improve on Construction~\ref{twofaceconst}.  However, this difference is something intrinsic to flag complexes.

\begin{lemma}  \label{flagkklemma}
Given positive integers $m$ and $k$, there is a unique way to pick integers $s \geq 0$ and $n_k, n_{k-1}, a_{k-1}, a_{k-2} \dots, a_{k-s}$ such that
$$m = {n_k \choose k} + {n_{k-1} \choose k-1} + {a_{k-1} \choose k-1} + {a_{k-2} \choose k-2} \dots + {a_{k-s} \choose k-s}$$
and $n_k > n_{k-1}$, $a_{k-1} > a_{k-2} > \dots > n_{k-s} \geq k-s > 0$, and
$${n_{k-1} \choose k-2} > {a_{k-1} \choose k-1} + {a_{k-2} \choose k-2} \dots + {a_{k-s} \choose k-s}.$$
\end{lemma}

This lemma is very similar to Lemma~\ref{kklemma}, with the key difference being two consecutive $k-1$ terms here.  Indeed, this lemma gives exactly the same values of $n_k$ and $n_{k-1}$ as Lemma~\ref{kklemma}.  The lemma is only interesting in the case of $k \geq 3$, as otherwise, the $a$ terms vanish and we are left with Lemma~\ref{kklemma}.  A previous paper of this author gave some stronger bounds for flag complexes than those of the Kruskal-Katona theorem.  \cite{mysecond}

\begin{theorem}  \label{flagkktheorem}
Let $m > 0$ and $k > p \geq 1$ be integers.  Let
$$m = {n_k \choose k} + {n_{k-1} \choose k-1} + {a_{k-1} \choose k-1} + {a_{k-2} \choose k-2} + \dots + {a_{k-s} \choose k-s}$$
as in Lemma~\ref{flagkklemma}.  Let $r = n_k - 1$ and let
$$m = {b_k \choose k}_r + {b_{k-1} \choose k-1}_{r-1} + \dots + {b_{k-s} \choose k-s}_{r-s}$$
as in Lemma~\ref{colorcan}.  If $\Delta$ is a flag complex with $f_{k-1}(\Delta) = m$, then either
$$f_{p-1}(\Delta) \geq {n_k \choose p} + {n_{k-1} \choose p-1} + {a_{k-1} \choose p-1} + {a_{k-2} \choose p-2} \dots + {a_{k-s} \choose p-s}$$
or else
$$f_{p-1}(\Delta) \geq {b_k \choose p}_r + {b_{k-1} \choose p-1}_{r-1} + \dots + {b_{k-s} \choose p-s}_{r-s}.$$
\end{theorem}

Where the theorem comes from is that if $\Delta$ has a simplex on $n_k$ vertices, then $f_{p-1}(\Delta)$ must satisfy the first bound, and if not, then it must satisfy the second bound by Theorem~\ref{flagcolored}.  In particular, $f_{p-1}(\Delta)$ must be larger then whichever is the smaller of the two prospective lower bounds.

The first lower bound is often the smaller of the two, and in particular, if $n_{k-1}$ is small relative to $n_k$, then the first lower bound is always the smaller of the two.  As $m \to \infty$, this still allows $a_{k-1}$ to become arbitrarily large while the first bound of Theorem~\ref{flagkktheorem} is still the smaller of the two; it simply has to grow much more slowly than $n_k$.

What makes this theorem important for our purposes is that if we let $q = {a_{k-1} \choose k-1} + {a_{k-2} \choose k-2} \dots + {a_{k-s} \choose k-s}$, then the first lower bound is $$f_{p-1}(\Delta) \geq {n_k \choose p} + {n_{k-1} \choose p-1} + g_{p-1}^{k-1}(q).$$
From Lemma~\ref{kklimit}, we get
$$\lim_{q \to \infty} {g_{p-1}^{k-1}(q) \over q^{{p-1 \over k-1}}} = {((k-1)!)^{{p-1 \over k-1}} \over (p-1)!}.$$
Hence, this is the best possible limit for the third-order terms that we could hope for in a construction that deals with flag complexes.  This is exactly the same limit as we found in Theorem~\ref{twofaceprop}.  An improvement to Construction~\ref{twofaceconst} will either leave this third-order limit alone or else make the limit no longer exist. Furthermore, even if an improvement did make the limit not exist, the limit inferior would still be this value.

Theorem~\ref{flagkktheorem} also suggests an improvement to Constructions~\ref{twofaceconst} and~\ref{dimconst}.  If $s=1$ in Lemma~\ref{flagkklemma}, then Construction~\ref{twofaceconst} can attain the first bound of Theorem~\ref{flagkktheorem}.  This bound is often attained by a flag complex if $s=2$, as well.  In Construction~\ref{twofaceconst}, we started with a face on $n_k$ vertices, and then added two extra vertices that were adjacent to $n_{k-1}$ and $a_{k-1}$ of the original $n_k$ vertices, respectively.  If we make these two extra vertices adjacent to each other and choose the vertices of the original $n_k$ to which they are adjacent such that they have $a_{k-2}$ common neighbors, then this adds ${a_{k-2} \choose k-2}$ faces of dimension $k-1$ and ${a_{k-2} \choose p-2}$ faces of dimension $p-1$, just as in the first bound of Theorem~\ref{flagkktheorem}.

This cannot always be done, however.  It requires $a_{k-2}$ vertices adjacent to both of the extra vertices, $n_{k-1} - a_{k-2}$ vertices adjacent to the first extra vertex but not the second, and $a_{k-1} - a_{k-2}$ vertices adjacent to the second extra vertex but not the first.  This requires $n_{k-1} + a_{k-1} - a_{k-2}$ vertices among the first $n_k$, which is only possible if $n_k \geq n_{k-1} + a_{k-1} - a_{k-2}$.  If this inequality holds, then we can make this adjustment to the construction.

This tweak is usually an improvement, too.  The big difference between the bound of the Kruskal-Katona theorem and what we can construct in Construction~\ref{twofaceconst} is that the latter has to keep adding terms of the form ${c \choose k-1}$, while the former can put smaller values rather than $k-1$ for all terms past the second.  This allows the same number of faces of dimension $k-1$ with fewer faces of smaller dimensions.  Having a $k-2$ rather than $k-1$ for a lot of the terms usually allows fewer faces as well.

It is fairly obvious how to make the same tweak to Construction~\ref{dimconst}.  We incorporate this adjustment in Construction~\ref{mainconst} in the next section.  We did not do so in Constructions~\ref{twofaceconst} and~\ref{dimconst} because it would not change the limits of Theorems~\ref{twofaceprop} and~\ref{dimprop}, but it would make the analysis much messier.

We would also like to do a quick analysis on the limit of Theorem~\ref{dimprop}.

\begin{definition}
Let $m > 0$ and $r \geq k > p \geq 1$ be integers.  Let
$$m = {n_k \choose k}_r + {n_{k-1} \choose k-1}_{r-1} + \dots + {n_{k-s} \choose k-s}_{r-s}$$
as in Lemma~\ref{colorcan}.  Define
$$j_{p,r}^k(m) = {n_k \choose p}_r + {n_{k-1} \choose p-1}_{r-1} + \dots + {n_{k-s} \choose p-s}_{r-s}.$$
\end{definition}

With this definition, Theorem~\ref{coloredkk} essentially states that if $\Delta$ is an $r$-colored simplicial complex with $f_{k-1}(\Delta) = m$, then $f_{p-1}(\Delta) \geq j_{p,r}^k(m)$.

\begin{lemma} \label{colorkklimit}
Let $r \geq k > p \geq 1$ be integers and let $f(x)$ be a function whose domain is the natural numbers such that for positive integers $m$ and $n$, if ${nr \choose k}_r \leq m \leq {(n+1)r \choose k}_r$, then ${nr \choose p}_r \leq f(m) \leq {(n+1)r \choose p}_r$.  Then
$$\lim_{m \to \infty} {f(m) \over m^{p \over k}} = {r \choose p}{r \choose k}^{-{p \over k}}.$$
\end{lemma}

As with Lemma~\ref{kklimit}, the use of this lemma is that $j_{p,r}^k(m)$ is a function that satisfies the requirements.  The proof here is analogous to the proof of that lemma.

\proof  For each integer $m$, let $n_r^m$ be the unique integer such that ${n_r^m r \choose k}_r \leq m < {(n_r^m +1)r \choose k}_r$.  We can compute
\begin{eqnarray*}
{{n_r^m r \choose p}_r \over {(n_r^m +1)r \choose k}_r^{{p \over k}}} \leq & {f(m) \over m^{{p \over k}}} & \leq {{(n_r^m +1) r \choose p}_r \over {n_r^m r \choose k}_r^{{p \over k}}} \\ {{r \choose p}(n_r^m)^p \over \big({r \choose k} (n_r^m +1)^k\big)^{{p \over k}}} \leq & {f(m) \over m^{{p \over k}}} & \leq {{r \choose p}(n_r^m + 1)^p \over \big({r \choose k} (n_r^m)^k\big)^{{p \over k}}} \\ {r \choose p}{r \choose k}^{-{p \over k}} \bigg({n_r^m \over n_r^m + 1}\bigg)^p \leq & {f(m) \over m^{{p \over k}}} & \leq {r \choose p}{r \choose k}^{-{p \over k}} \bigg({n_r^m + 1 \over n_r^m}\bigg)^p.
\end{eqnarray*}
As $m \to \infty$, we get $n_r^m \to \infty$ also.  As $n_r^m \to \infty$, we have both ${n_r^m \over n_r^m + 1} \to 1$ and ${n_r^m + 1 \over n_r^m} \to 1$.  Thus, both ends of the last big inequality go to ${r \choose p}{r \choose k}^{-{p \over k}}$ as $m \to \infty$.  Therefore, so does ${f(m) \over m^{{p \over k}}}$.  \endproof

It follows immediately from this lemma that
$$\lim_{m \to \infty} {j_{p,r}^k(m) \over m^{{p \over k}}} = {r \choose p}{r \choose k}^{-{p \over k}}.$$
Likewise, it is not difficult to show that the analogous first-order limit for Construction~\ref{dimconst} is also ${r \choose p}{r \choose k}^{-{p \over k}}$.  Similarly, if we chop off the first term in both cases, the lemma gives
$$\lim_{m \to \infty} {j_{p,r}^k(m) - {n_k \choose k}_r \over m^{{p-1 \over k-1}}} = {r-1 \choose p-1}{r-1 \choose k-1}^{-{p-1 \over k-1}}.$$
It is not difficult to show that the analogous second-order limit of Construction~\ref{dimconst} is the same as this.

The limits only differ if we chop off the first two terms.  Here, we get
$$\lim_{m \to \infty} {j_{p,r}^k(m) - {n_k \choose k}_r - {n_{k-1} \choose k-1}_{r-1} \over m^{{p-2 \over k-2}}} = {r-2 \choose p-2}{r-2 \choose k-2}^{-{p-2 \over k-2}}.$$
This limit differs from that of Theorem~\ref{dimprop}, and as in the case with the dimension not restricted, the key difference is not the different limit value, but rather the different power of $m$ in the denominator of the quantity of which we take the limit.

As with the situation where we did not have any restrictions on the dimension, it is likely that this difference is due to properties intrinsic to flag complexes, rather than a result of the construction not being very good.  In this scenario, however, we do not have a theorem analogous to Theorem~\ref{flagkktheorem}.

\section{The main construction}

If one has a proposed face vector for a flag complex and wants to find a flag complex with the given face vector, it probably isn't a case of wanting to match two particular face numbers.  More commonly, one would want a flag complex that has the entire desired face vector.  In this section, we discuss how to construct a flag complex that does exactly that, and borrow heavily from the approaches of earlier sections.

The intuitive idea of the construction is that we use the highest dimensional faces first.  In Construction~\ref{dimconst}, we tried to use exactly the prescribed number of $(k-1)$-dimensional faces with as few $(p-1)$-dimensional faces as possible.  The construction depended only on $k$, however, and not on $p$, and did pretty well at using relatively few faces of all smaller dimensions.

If we have a proposed face vector for a complex of dimension $d-1$, we essentially start with Construction~\ref{dimconst}, adjusted as discussed in Section~4, using $r = k = d$ to add in all of the $(d-1)$-dimensional faces while using few faces of smaller dimensions.  Next, we add however many more $(d-2)$-dimensional faces are needed, essentially using Construction~\ref{dimconst} again, this time using $r = k = d-1$ to only add more faces of dimension up to $d-2$, while again using as few faces of smaller dimensions as we can.  We repeat this for each dimension until all of the faces have been added.

\begin{construction}  \label{mainconst}
\textup{Let $(1, c_1, c_2, \dots, c_d)$ be a vector of positive integers.  We wish to create a flag complex $\Delta$ of dimension $d-1$ with $f_{i-1}(\Delta) = c_i$ for all $1 \leq i \leq d$.  Define $n_0^d$ to be the unique integer such that ${n_0^d \choose d}_d \leq c_d < {n_0^d + 1 \choose d}_d$.  Define $m_0^d = c_d - {n_0^d \choose k}_d$.}

\textup{Define $m_i^d$, $n_i^d$, $p_i^d$, and $q_i^d$ for $i \geq 1$ recursively as follows.  If $m_{i-1}^d = 0$, then let $z^d = i-1$ and stop.  If $m_{i-1}^d > 0$, then let $n_i^d$ be the unique integer such that ${n_i^d \choose d-1}_{d-1} \leq m_{i-1}^d < {n_i^d + 1 \choose d-1}_{d-1}$.  Let $q_i^d$ be the unique integer such that ${q_i^d \choose d-2}_{d-2} \leq m_{i-1}^d - {n_i^d \choose d-1}_{d-1} < {q_i^d + 1 \choose d-2}_{d-2}$.  If
$$n_0^d - \bigg\lfloor {n_0^d \over d} \bigg\rfloor - \bigg\lfloor {n_0^d + 1 \over d} \bigg\rfloor + q_i^d \geq n_{i-1}^d - \bigg\lfloor {n_{i-1}^d \over d-1} \bigg\rfloor + n_i^d - \bigg\lfloor {n_i^d \over d-1} \bigg\rfloor,$$
then let $p_i^d = q_i^d$.  If $n_0^d + q_i^d < n_{i-1}^d + n_i^d$, then let $p_i^d = -1$.  Let $m_i^d = m_{i-1}^d - {n_i^d \choose d-1}_{d-1} - {p_i^d \choose d-2}_{d-2}$.}

\textup{Construct a complex $\Gamma^d$ as the clique complex of a graph that starts with the Tur\'{a}n graph $T_{n_0,d}$, and then adds $z$ additional vertices, $v_1, v_2, \dots, v_z$, such that the link of $v_i$ is a Tur\'{a}n graph $T_{n_i,d-1}$ that is a subgraph of the original Tur\'{a}n graph $T_{n_0,d}$ plus, if $p_i > 0$, the vertex $v_{i-1}$ adjacent to a Tur\'{a}n graph $T_{p_i, d-2}$ that is a subgraph of the $T_{n_i,d-1}$ in the link.}

\textup{It is easy to compute that for each $k$,}
$$f_{k-1}(\Gamma^d) = {n_0 \choose k}_d + {n_1 \choose k-1}_{d-1} + {p_1 \choose k-2}_{d-2} + \dots + {n_{z^d} \choose k-1}_{d-1} + {p_{z^d} \choose k-2}_{d-2}$$
\textup{It is also easy to see that}
$$m_{i-1}^d = {n_i \choose d-1}_{d-1} + {p_i \choose d-2}_{d-2} + \dots + {n_{z^d} \choose d-1}_{d-1} + {p_{z^d} \choose d-2}_{d-2},$$
\textup{by starting with $i = z$ and working backwards.  From this, it follows that $f_{d-1}(\Gamma^d) = m$.}

\textup{Now we define $\Delta^d = \Gamma^d$ and repeat the basic process to create a complex $\Gamma^{d-1}$.  If $f_{d-2}(\Delta^d) > c_{d-1}$, then we have already used more faces of dimension $d-2$ than allowed, so the construction fails.  Otherwise, we let $n_0^{d-1}$ be the unique integer such that}
$${n_0^{d-1} \choose d-1}_{d-1} \leq c_{d-1} - f_{d-2}(\Delta^d) + {n_0^d - \big\lfloor {n_0^d \over d} \big\rfloor \choose d-1}_{d-1} < {n_0^{d-1} + 1 \choose d-1}_{d-1}.$$
\textup{Also, we define}
$$m_0^{d-1} = c_{d-1} - f_{d-2}(\Delta^d) + {n_0^d - \big\lfloor {n_0^d \over d} \big\rfloor \choose d-1}_{d-1} - {n_0^{d-1} \choose d-1}_{d-1}.$$

\textup{Next, we define $m_i^{d-1}$, $n_i^{d-1}$, $p_i^{d-1}$, and $q_i^{d-1}$ for $i \geq 1$ recursively in the same manner as we defined $m_i^d$, $n_i^d$, $p_i^d$, and $q_i^d$ in the second paragraph of this construction, except using $d-1$ everywhere rather than $d$.  Likewise, we construct a flag complex $\Gamma^{d-1}$ just like in the third paragraph of this construction.  We can compute $f_{k-1}(\Gamma^{d-1})$ just as we did $f_{k-1}(\Gamma^d)$.}

\textup{Now we define $\Delta^{d-1} = \Gamma^{d-1} \cup \Delta^d$, where the intersection of the two complexes is the $T_{n_0^d - \lfloor{n_0^d \over d}\rfloor,d-1}$ that uses all vertices of the $T_{n_0^d,d}$ at the start of building $\Gamma^d$ of the first $d-1$ colors, as well as the first several vertices of each color of a $T_{n_0^d - \lfloor{n_0^d \over d}\rfloor,d-1}$ subgraph of the $T_{n_0^{d-1},d-1}$ at the start of building $\Gamma^{d-1}$.  From this, we can compute by an easy inclusion-exclusion argument that}
$$f_{k-1}(\Delta^{d-1}) = f_{k-1}(\Gamma^{d-1}) + f_{k-1}(\Delta^d) - {n_0^d - \lfloor{n_0^d \over d}\rfloor \choose d-1}_d.$$

\textup{Because $\Gamma^{d-1}$ does not contain any faces of dimension $d-1$, we get $f_{d-1}(\Delta^{d-1}) = f_{d-1}(\Delta^d) = c_d$.  $\Gamma^{d-1}$ is defined to add just enough faces of dimension $d-1$ to make $f_{d-2}(\Delta^{d-1}) = c_{d-1}$.  Thus, $\Delta^{d-1}$ has exactly the right number of faces for the last two numbers in its face vector.  We repeat this, constructing and tacking on $\Gamma^{d-2}$, $\Gamma^{d-3}$, and so forth, and each new $\Delta^i$ has the appropriate face number for one more number at the end of its vector.  We stop when either the construction fails because we have already used too many faces of some dimension $i-1$ in $\Delta^{i+1}$ or when we construct $\Delta^1$, which will have exactly the desired face vector.  \endproof}
\end{construction}

This construction works for most of the same reasons as Constructions~\ref{twofaceconst} and~\ref{dimconst}.  There are some additional comments that are necessary, however.  First, throughout the construction, the superscript $d$'s are indices, not exponents.  Next, we use the convention ${n \choose k}_r = 0$ if $n < 0$ or $k < 0$.  Setting $p_i^j = -1$ is simply a way to say that a new vertex shouldn't be adjacent to the previous vertex at all.  Also, that ${n \choose k}_r = 0$ if $r < k$ follows from the definition.

The larger potential issue that needs to be addressed is that at certain places, the construction requires that two particular complexes have as their intersection a particular Tur\'{a}n graph.  We need to ensure that it is always possible to do this.

In the third paragraph of the construction, if we take a $T_{n_0,d}$, we can create a $T_{n_0+1,d}$ by adding another vertex whose link is a $T_{n_0 - \lfloor{n_0^d \over d}\rfloor,d-1}$.  Thus, there is room to add a vertex whose link is any subgraph of a $T_{n_0 - \lfloor{n_0^d \over d}\rfloor,d-1}$.  We must have $n_i < n_0 - \big\lfloor{n_0^d \over d}\big\rfloor$, for otherwise, $n_0$ would have been larger.  Thus, there is room for each vertex to have the appropriate link.

The bigger issue is whether the link of the edge $\{v_{i-1}, v_i\}$ can be the appropriate graph $T_{p_i, d-2}$.  Suppose that parts 1 through $w$ of the graph $T_{n_0,d}$ have  $\big\lceil{n_0^d \over d}\big\rceil$ vertices, while parts $w+1$ through $d$ have $\big\lfloor{n_0^d \over d}\big\rfloor$.  When choosing the vertices of the $T_{n_0,d}$ for $v_1$ to be adjacent to, we can start with a vertex of color 1, then one of color 2, and so forth, until we have one of each color except for $d$, and then repeat the cycle.  Do this until all vertices for the link of $v_1$ are chosen.

When choosing the vertices of the $T_{n_0,d}$ for the link of $v_2$, add vertices in the same order as for the link of $v_1$ except that we skip the vertices of color $d-1$ until the $p_2^d$ vertices in the link of the edge $\{v_1, v_2\}$ have been chosen.  Once this is done, for the remaining $n_2^d - p_2^d$ vertices for the link of $v_2$, we wish to add $\big\lfloor {p_2^d \over d-1} \big\rfloor$ vertices of color $d$.  Note that we can do this, as if $n_2^d < p_2^d + \big\lfloor {p_2^d \over d-1} \big\rfloor$, then we would have chosen $n_2^d$ to be larger.

Finally, we add vertices cycling through all of the colors except for $d-1$, such that no color ever has more than one vertex more than any other color except for color $d-1$.  We choose the order of colors to add the remaining vertices to the link of $v_2$ as follows.  The first colors to get an extra vertex are colors 1 through $w$ with $\big\lfloor {p_2^d \over d-2} \big\rfloor$ vertices in the link of $\{v_1, v_2\}$.  The next colors to get an extra vertex are colors $w+1$ through $d-2$ or color $d$ with $\big\lfloor {p_2^d \over d-2} \big\rfloor$ vertices in the link of $\{v_1, v_2\}$.  After that, we use colors 1 through $w$ with $\big\lfloor {p_2^d \over d-2} \big\rfloor + 1$ vertices in the link of $\{v_1, v_2\}$.  Finally, we use the remaining colors.  At least one of these sets of colors will be empty; when we hit an empty set of colors, we just skip it and move on.

We can easily check that this order of adding the remaining vertices of various colors will use all vertices available to be added, except those of color $d-1$, before we try to add another vertex of some color for which all vertices not adjacent to $v_1$ are already adjacent to $v_2$.

The construction makes $v_1$ is adjacent to $\big\lfloor {n_{i-1}^d \over d-1} \big\rfloor$ vertices of color $d-1$ and none of color $d$, so it is adjacent to $n_{i-1}^d - \big\lfloor {n_{i-1}^d \over d-1} \big\rfloor$ vertices among the first $d-2$ colors.  Similarly, $v_2$ is adjacent to $\big\lfloor {n_i^d \over d-1} \big\rfloor$ vertices of color $d$ and none of color $d-1$, so it is adjacent to $n_i^d - \big\lfloor {n_i^d \over d-1} \big\rfloor$ vertices among the first $d-2$ colors.  Summing these two figures gives the number of vertices adjacent to $v_1$ but not $v_2$ plus the number adjacent to $v_2$ but not $v_1$ plus the number of vertices adjacent to both.  There are $n_0^d - \big\lfloor {n_0^d \over d} \big\rfloor - \big\lfloor {n_0^d + 1 \over d} \big\rfloor$ vertices among the first two colors in all, and we want exactly $p_2^d$ of them to be adjacent to both $v_1$ and $v_2$.  Since we must have
$$n_0^d - \bigg\lfloor {n_0^d \over d} \bigg\rfloor - \bigg\lfloor {n_0^d + 1 \over d} \bigg\rfloor + p_i^d \geq n_{i-1}^d - \bigg\lfloor {n_{i-1}^d \over d-1} \bigg\rfloor + n_i^d - \bigg\lfloor {n_i^d \over d-1} \bigg\rfloor$$
in order for $v_1$ and $v_2$ to be adjacent to each other at all, it is possible to choose vertices such that there are exactly $p_2^d$ vertices adjacent to both $v_1$ and $v_2$.

It is easy to check that the order of adding the remaining vertices of various colors ensures after each step, $v_2$ is adjacent to at least as many vertices of each color with $\lceil{n_0^d \over d}\rfloor$ vertices in the original $T_{n_0,d}$ as of each color with $\lfloor{n_0^d \over d}\rfloor$ vertices in the original $T_{n_0,d}$, except for color $d-1$, but never by more than one extra vertex.  Thus, when we go to add $v_3$, we can renumber the colors of the original $T_{n_0,d}$ such that $v_2$ satisfies the same condition on colors as $v_1$.

The construction works for constructing $\Gamma^{d-1}$, $\Gamma^{d-2}$, and so on by the same argument as for $\Gamma^d$.

The other thing that we need to check is that when we try to form $\Delta^i$ by intersecting $\Gamma^i$ with $\Delta^{i+1}$, both $\Gamma^i$ and $\Delta^{i+1}$ contain the appropriate graph $T_{n_0^{i+1} - \lfloor{n_0^{i+1} \over i+1}\rfloor,i}$ that we can use as their intersection.  Here, we note that $\Delta_{i+1}$ contains a $T_{n_0^{i+1},i+1}$, and removing the vertices of the last color yields a $T_{n_0^{i+1} - \big\lfloor {n_0^{i+1} \over i+1} \big\rfloor,i}$.  Furthermore, $\Gamma^i$ contains a $T_{n_0^i,i}$, and from the choice of $n_0^i$, we must have $n_0^i \geq n_0^{i+1} - \Big\lfloor {n_0^{i+1} \over i+1} \Big\rfloor$ or else the construction would have failed when we tried to add the faces of dimension $i-1$.  We can take the first $n_0^{i+1} - \Big\lfloor {n_0^{i+1} \over i+1} \Big\rfloor$ vertices from this graph to get the desired $T_{n_0^{i+1} - \big\lfloor {n_0^{i+1} \over i+1} \big\rfloor,i}$ subgraph.

Construction~\ref{mainconst} is rather complicated, so perhaps an example will illustrate how it works.

\begin{example}  \label{genexample}
\textup{We wish to construct a flag complex with face vector $(1, 100, 1000, 2000)$ by the method of Construction~\ref{mainconst}.  This complex will have dimension two, so $d = 3$.  We can compute ${37 \choose 3}_3 = 1872$ and ${38 \choose 3}_3 = 2028$, so $n_0^3 = 37$ and we start with a graph $T_{37, 3}$.  This gives $m_0^3 = 2000 - 1872 = 128$.  We compute ${22 \choose 2}_2 = 121$ and ${23 \choose 2}_2 = 132$, so $n_1^3 = 22$.  We add an extra vertex $v_1$ and make it adjacent to the first 11 vertices of each of the first two colors of the $T_{37, 3}$.  Since $p_1^j = -1$ always, we can skip computing $q_1^3$.  Next, we compute $m_1^3 = 128 - 121 = 7$.  We compute ${6 \choose 2}_2 = 9$ and ${5 \choose 2}_2 = 6$, so $n_2^3 = 5$.  We next compute ${1 \choose 1}_1 = 1$, so $q_2^3 = 1$.  To see whether $p_2^3$ is $-1$ or $q_2^3$, we check whether}
$$37 - \bigg\lfloor {37 \over 3} \bigg\rfloor - \bigg\lfloor {37 + 1 \over 3} \bigg\rfloor + 1 \geq 22 - \bigg\lfloor {22 \over 2} \bigg\rfloor + 5 - \bigg\lfloor {5 \over 2} \bigg\rfloor.$$
\textup{This simplifies to $14 \geq 14$, which is true.  Hence, $p_2^3 = 1$.  We make $v_2$ adjacent to the first two vertices of color 3 and the last three vertices of color 1.  This makes the eleventh vertex of color 1 in the $T_{37, 3}$ adjacent to both $v_1$ and $v_2$.  We also make $v_1$ and $v_2$ adjacent to each other.  At this point, we compute $m_2^3 = 0$, so $z^3 = 2$ and we are done with $\Gamma^3$ and hence $\Delta^3$.  At this point, we may wish to stop to compute}
\begin{eqnarray*}
f_2(\Delta^3) & = & {37 \choose 3}_3 + {22 \choose 2}_2 + {5 \choose 2}_2 + {1 \choose 1}_1 = 2000 \\ f_1(\Delta^3) & = & {37 \choose 2}_3 + {22 \choose 1}_2 + {5 \choose 1}_2 + {1 \choose 0}_1 = 484 \\ f_0(\Delta^3) & = & {37 \choose 1}_3 + {22 \choose 0}_2 + {5 \choose 0}_2 + {1 \choose -1}_1 = 39.
\end{eqnarray*}

\textup{We have $f_1(\Delta^3) \leq 1000$ and $f_0(\Delta^3) \leq 100$, so the construction has not yet failed.  If we restrict the graph $T_{37, 3}$ to the vertices of the first two colors, it is a $T_{25, 2}$.  This has ${25 \choose 2}_2 = 156$ edges, so we want for $\Gamma^2$ to contain $1000 - 484 + 156 = 672$ edges.  Toward this end, we compute ${51 \choose 2}_2 = 650$ and ${52 \choose 2}_2 = 676$, so $n_0^2 = 51$.  We start building $\Gamma^2$ with a graph $T_{51,2}$.  This gives us $m_0^2 = 672 - 650 = 22$.  We compute ${22 \choose 1}_1 = 22$, so we can finish $\Gamma^2$ by adding another vertex and making it adjacent to 22 vertices of color 1.  We compute $f_1(\Gamma^2) = 672$ and $f_0(\Gamma^2) = 52$.}

\textup{The intersection of $\Gamma^2$ and $\Delta^3$ in $\Delta^2$ is $T_{25, 2}$.  Thus, we can compute}
\begin{eqnarray*}
f_1(\Delta^2) & = & f_1(\Delta^3) + f_1(\Gamma^2) - {25 \choose 2}_2 = 484 + 672 - 156 = 1000 \\ f_0(\Delta^2) & = & f_0(\Delta^3) + f_0(\Gamma^2) - {25 \choose 1}_2 = 39 + 52 - 25 = 66
\end{eqnarray*}

\textup{From here, we only need to add more vertices.  It is easy to form $\Delta^1$ by adding 34 isolated vertices to $\Delta^2$.  The face vector of $\Delta^1$ is $(1, 100, 1000, 2000)$, as desired.}
\end{example}

Another point to be made is that, while less efficient, it often works just fine to skip the step of making extra vertices adjacent to each other.  This makes the construction simpler, though it also makes it no longer work in some borderline cases.

While Construction~\ref{mainconst} often works, there is one important adjustment that sometimes needs to be made.  If the complex is of dimension $d-1$ and the main obstruction to making it work is using enough faces of dimension $k-1$ without using too many faces of dimension $p-1$ for some $d > k > p$, it can be beneficial, rather than starting with a complete $d$-partite graph with the vertices divided among the parts as evenly as possible, to instead put more vertices into the first $k$ parts and fewer vertices into the last $d-k$ parts.  While this typically isn't the most efficient way to get fewer faces of dimensions $k-1$ and $p-1$ for a given number of faces of dimension $d-1$, that may not matter if there aren't that many faces of dimension $d-1$.  This tweak can sometimes allow fewer faces of dimension $p-1$ for a given number of faces of dimension $k-1$.  Perhaps an example will illustrate the point.

\begin{example}  \label{unbalancedexample}
\textup{We wish to construct a flag complex with face vector $(1, 62, 1161, 5832)$.  We can compute ${54 \choose 3}_3 = 5832$, so $\Delta^3 = \Gamma^3 = T_{54, 3}$.  We compute $f_1(\Delta^3) = 972$ and $f_0(\Delta^3) = 54$.  If we remove color 3 from $\Delta^3$, we are left with a $T_{36,2}$.  We thus want $f_1(\Gamma^2) = 1161 - 972 + {36 \choose 2}_2 = 513$.  We compute ${45 \choose 2}_2 = 506$ and ${46 \choose 2}_2 = 529$, so $n_0^2 = 45$.  This gives us $m_0^2 = 7$.  We compute ${7 \choose 1}_1 = 7$, so $n_1^2 = 7$.  We form $\Gamma^2$ by starting with a $T_{45, 2}$ and then adding one extra vertex adjacent to seven vertices of the first color.  This gives us $f_2(\Delta^2) = 5832$ and $f_1(\Delta^2) = 1161$.  However,}
$$f_0(\Delta^2) = f_0(\Delta^3) + f_0(\Gamma^2) - {36 \choose 1}_2 = 54 + 46 - 36 = 64.$$
\textup{We are only allowed 62 vertices, and have already used 64, so the construction fails.}

\textup{There is a complex with the desired face vector, however.  The clique complex of a complete tripartite graph with 27 vertices on each of the first two parts and 8 vertices on the third part has exactly the desired face vector.  If we started with this graph as $\Delta^3$, the construction would work.}
\end{example}

It is important to observe that this is not a case of discreteness getting in the way.  There was only one vertex added to a $\Gamma^i$ apart from the initial Tur\'{a}n graphs, and we still used two vertices too many.  This is a case where allocating the vertices into parts unevenly is essential.

This naturally raises the question of when one should allocate vertices into parts unevenly, and how unevenly to do so.  Both of these issues are addressed in the next section.

Once the uneven allocation of vertices into parts is done, it is fairly clear how to construct a $\Gamma^i$ from this as in Construction~\ref{mainconst}.  The formula to determine when two extra vertices can be adjacent to each other breaks down, but given a particular desired face vector, an ad hoc approach to see if it's possible to make two extra vertices adjacent to each other works just fine.  Other than the first two extra vertices, it usually is pretty easy to make subsequent pairs adjacent to each other.

\section{Some conjectures}

In this section, we present two conjectures that, together, would constitute great progress toward characterizing the face vectors of flag complexes, and argue as to why they are likely to be true.  These conjectures rely on the existence of some constants, so we subsequently discuss how to find the constants if they exist.  We then return to the question of when to put more vertices in some parts than others in applying Construction~\ref{mainconst}.

\begin{conjecture}  \label{bigconj}
Let $(1, c_1, c_2, \dots, c_d)$ be the face vector of a flag complex.  Then there are nonnegative real numbers $a_j^k$ for all integers $1 \leq j \leq k \leq d$ such that
\begin{enumerate}
\item for every $1 \leq j \leq d$,
$$c_j = \sum_{i_1 < i_2 < \dots < i_j} a_{i_1}^{i_j} a_{i_2}^{i_j} \dots a_{i_j}^{i_j},$$
\item for every $i$, $j$, and $k$ with $i < j$, $a_i^k \geq a_j^k$, and
\item for every $i$, $j$, and $k$ with $i < k$, $a_j^i \geq a_j^k$.
\end{enumerate}
\end{conjecture}

The intuition here is that we drop the requirement that vertices be discrete, and instead allow things like a complete tripartite graph with 9.5, ${25 \over 3}$, and 4.7 vertices in the respective parts.  The edges connecting vertices of color $j$ to those of color $k$ connect the first $a_j^k$ vertices of color $j$ to the first $a_k^k$ vertices of color $k$.  The first condition of the lemma says that this produces the appropriate face vector, as $a_{i_1}^{i_j} a_{i_2}^{i_j} \dots a_{i_j}^{i_j}$ is the number of faces of dimension $j-1$ using a vertex of each of colors $i_1, i_2, \dots, i_j$.  The sum is over all ways to pick the $j$ colors.  The second and third conditions are essentially shifting conditions.  The second condition basically says that we can decide that the color that gets used the most is color 1, the color that gets used the second most is color 2, and so forth.

The real intuition of this lemma is that if we allow non-integer numbers of vertices in each part, the converse would be true, essentially by Construction~\ref{mainconst}, adjusted to allow uneven numbers of vertices in the various parts.  Following the approach of Construction~\ref{mainconst}, we start with a complete $d$-partite graph with $a_1^d$ vertices of color 1, $a_2^d$ vertices of color 2, and so forth. Next, we add a complete $(d-1)$-partite graph with $a_1^{d-1}$ vertices of color 1, $a_2^{d-1}$ vertices of color 2, and so forth.  The third condition of the lemma says that this includes the complete $(d-1)$-partite graph that we would get by removing the vertices of color $d$ from the original $d$-partite graph.  We then add a complete $(d-2)$-partite graph, and so forth.

\begin{conjecture}  \label{converseconj}
Let $(1, c_1, c_2, \dots, c_d)$ be a vector of positive integers and let $a_j^k$ be real numbers for all integers $1 \leq j \leq k \leq d$ such that
\begin{enumerate}
\item for every $1 \leq j \leq d$,
$$c_j = \sum_{i_1 < i_2 < \dots < i_j} a_{i_1}^{i_j} a_{i_2}^{i_j} \dots a_{i_j}^{i_j},$$
\item for every $i$, $j$, and $k$ with $i < j$, $a_i^k \geq a_j^k$, and
\item for every $i$, $j$, and $k$ with $i < k$, $a_j^i > a_j^k$.
\end{enumerate}
Then there is an integer $q$ such that $(1, q c_1, q^2 c_2, \dots, q^d c_d)$ is the face vector of a flag complex.
\end{conjecture}

Furthermore, this conjecture is true.  If all of the $a$'s are rational, then we can pick $q$ to be the least common multiple of the denominators, so that $qa_j^k$ is an integer for all choices of $j$ and $k$.  We can then apply essentially Construction~\ref{mainconst}, except with uneven parts.  We start with a complete $d$-partite graph with $a_i^d$ vertices of color $i$.  Then we add a complete $(d-1)$-partite graph with $a_i^{d-1}$ vertices of color $i$, and so forth.  At the end of the construction, we will have a flag complex whose face vector is precisely $(1, q c_1, q^2 c_2, \dots, q^d c_d)$.

It is intentional that the third condition is a strict inequality in Conjecture~\ref{converseconj} and a weak inequality in Conjecture~\ref{bigconj}.  The strict inequality of Conjecture~\ref{converseconj} is necessary to allow some margin to account for discreteness doing weird things.  Since there are and ${d \choose 2}$ $a$'s and only $d$ $c$'s, if the $a$'s are real, then we should be able to perturb them a bit to make them rational while maintaining the various conditions of the conjecture.

An argument along the lines of Theorems~\ref{twofaceprop} and~\ref{dimprop} would limit how much discreteness could interfere, and likely prove Conjecture~\ref{converseconj}.  Furthermore, such a proof would probably give values of $q$ that will suffice, and if $q = 1$ works, then it would show that $(1, c_1, c_2, \dots, c_d)$ is the face vector of a flag complex.   Unfortunately, the details are likely to be long and messy, which is why Conjecture~\ref{converseconj} is only presented as a conjecture and not a theorem.

At the end of Section~5, we promised to return to the question of when it was best to start with a Tur\'{a}n graph in constructing some $\Gamma^i$ and when it was better to start with a graph that put significantly more vertices in some parts than others.  The answer to this is that finding how uneven to make the parts in the construction is very similar to finding the constants $a_j^k$ as in Conjecture~\ref{bigconj}.  When constructing $\Gamma^k$ as in Construction~\ref{mainconst}, one wants to start with about $a_j^k$ vertices of color $k$.

There could, of course, be many possible choices of $a_j^k$, and many ways to pick integers near them.  It helps to make the third condition of the conjecture far from equality, that is, $a_j^i - a_j^k$ should be ``large" for all relevant choices of $i$, $j$, and $k$.

This doesn't always work perfectly, of course, and one may need to tweak the construction a bit at times.  Furthermore, there are integer vectors that do not correspond to any flag complex, even though they satisfy the conditions of Conjecture~\ref{bigconj}.  It also isn't a deterministic method, but only says, here is roughly where to look and you might find something that works.

Even so, this only moves the problem from the discrete situation to the continuous one, and does not explain where the constants $a_j^k$ come from.  The simplest thing to try is to start with $a_1^d = a_2^d = \dots = a_d^d = (c_d)^{1 \over d}$. Next, we set $a_1^{d-1} = a_2^{d-1} = \dots = a_{d-1}^{d-1}$, and set the common value to whatever is needed to make the first condition of Conjecture~\ref{bigconj} work for $j = d-1$. We repeat this for $j = d-2$, then $j = d-3$, and so forth.  The second condition of the conjecture is trivially satisfied, and if the third is also, then we have our constants.

This approach is akin to trying Construction~\ref{mainconst} unmodified.  Sometimes it does work, as in Example~\ref{genexample}.  But sometimes it fails, even when an unbalanced allocation of vertices among the colors would succeed as in Example~\ref{unbalancedexample}.

The way that it can fail is when the simple approach above gives $a_1^i < a_1^{i+1}$, in contradiction of the third requirement of the conjecture.  When this happens, we know which particular dimension of faces we run out of prematurely.  What needs to be tweaked is to allow the same number of faces of all higher dimensions, while using up fewer faces of the dimension that we run out of early.  For this, we have a theorem, a proposition, and a conjecture.

\begin{theorem}  \label{equalvertices}
Let $\Delta$ be a flag complex of dimension $d-1$.  Let $p$ and $k$ be integers such that $1 \leq p < k \leq d$.  Then
$$f_{p-1}(\Delta) \geq {d \choose p} {d \choose k}^{-{p \over k}} (f_{k-1}(\Delta))^{{p \over k}}.$$
\end{theorem}

This theorem is basically a numerical approximation to Theorem~\ref{coloredkk}.  It is equivalent to \cite[Theorem 5.1]{balanced}, though formulated differently.

The intuition is that if we ignore discreteness, then one way to construct a $d$-colored flag complex $\Delta$ with a desired value of $f_{k-1}(\Delta)$ is to have a complete $d$-partite graph with $w$ vertices in each color.  One can readily compute that $f_{k-1}(\Delta) = w^k {d \choose k}$, from which $w = (f_{k-1}(\Delta))^{-{1 \over k}} {d \choose k}^{-{1 \over k}}$.  We can likewise compute
$$f_{p-1}(\Delta) = w^p {d \choose p} = {d \choose p} {d \choose k}^{-{p \over k}} (f_{k-1}(\Delta))^{{p \over k}}.$$
The theorem says that this is a lower bound for the number of faces of dimension $p-1$ that the complex must contain.

A bit of straightforward but messy analysis can convert this theorem to fit our setup that ignores discreteness.

\begin{proposition}  \label{contequalvert}
Let $1 \leq p < k \leq d$ be integers and let $m > 0$.  Let $a_1, a_2, \dots, a_d \geq 0$ be real numbers such that
$$m = \sum_{i_1 < i_2 < \dots < i_k} a_{i_1} a_{i_2} \dots a_{i_k}.$$
Then
$$\sum_{i_1 < i_2 < \dots < i_p} a_{i_1} a_{i_2} \dots a_{i_p} \geq {d \choose p} {d \choose k}^{-{p \over k}} m^{{p \over k}}.$$
\end{proposition}

\proof  Suppose not.  Then there must be some choice $a_1, a_2, \dots, a_d \geq 0$ which is a counterexample.  Let
$$c = {d \choose p} {d \choose k}^{-{p \over k}} m^{{p \over k}} - \sum_{i_1 < i_2 < \dots < i_p} a_{i_1} a_{i_2} \dots a_{i_p}.$$
Let $w = {d \choose p}^{-1} c$.  We have $w > 0$ because the $a's$ are a counterexample to the corollary.  For any choice of $i_1 < i_2 < \dots < i_k$, let
$$f_{i_1, i_2, \dots, i_p}(x) = (a_{i_1}-x) (a_{i_2}-x) \dots (a_{i_p}-x).$$
It is clear that $f_{i_1, i_2, \dots, i_p}(x)$ is a polynomial in $x$, and hence continuous in $x$.  Therefore, there is some $\delta > 0$ such that if $|x| < \delta$, then $|f_{i_1, i_2, \dots, i_p}(x) - f_{i_1, i_2, \dots, i_k}(0)| < w$.

Let $\delta_0$ be the smallest such choice of $\delta$ over all ${d \choose p}$ choices of $i_1 < i_2 < \dots < i_p$.  Let $z = \lceil {1 \over \delta_0} \rceil$.  For every integer $i$ with $1 \leq i \leq d$, let $b_i = {1 \over z} \lceil a_i z \rceil$.  We can compute
\begin{eqnarray*}
& & \sum_{i_1 < i_2 < \dots < i_p} b_{i_1} b_{i_2} \dots b_{i_p} \\ & \leq & \sum_{i_1 < i_2 < \dots < i_p} (a_{i_1} + \delta_0) (a_{i_2} + \delta_0) \dots (a_{i_p} + \delta_0) \\ & < & \sum_{i_1 < i_2 < \dots < i_p} (a_{i_1} a_{i_2} \dots a_{i_p} + w) \\ & = & -c + \sum_{i_1 < i_2 < \dots < i_p} a_{i_1} a_{i_2} \dots a_{i_p} \\ & = & {d \choose p} {d \choose k}^{-{p \over k}} m^{{p \over k}}
\end{eqnarray*}

Let $\Delta$ be a complete $d$-partite graph with $zb_1, zb_2, \dots zb_d$ vertices in its respective parts.  Then
\begin{eqnarray*}
f_{p-1}(\Delta) & = & \sum_{i_1 < i_2 < \dots < i_p} (zb_{i_1}) (zb_{i_2}) \dots (zb_{i_p}) \\ & = & z^p \sum_{i_1 < i_2 < \dots < i_p} b_{i_1} b_{i_2} \dots b_{i_p} \\ & < & z^p {d \choose p} {d \choose k}^{-{p \over k}} m^{{p \over k}} \\ & = & {d \choose p} {d \choose k}^{-{p \over k}} (z^k m)^{{p \over k}}.
\end{eqnarray*}
However, we can compute
\begin{eqnarray*}
f_{k-1}(\Delta) & = & \sum_{i_1 < i_2 < \dots < i_k} (zb_{i_1}) (zb_{i_2}) \dots (zb_{i_k}) \\ & = & z^k \sum_{i_1 < i_2 < \dots < i_k} b_{i_1} b_{i_2} \dots b_{i_k} \\ & \geq & z^k \sum_{i_1 < i_2 < \dots < i_k} a_{i_1} a_{i_2} \dots a_{i_k} \\ & = & z^k m
\end{eqnarray*}
Therefore, we have
$$f_{p-1}(\Delta) < {d \choose p} {d \choose k}^{-{p \over k}} (z^k m)^{{p \over k}} \leq {d \choose p} {d \choose k}^{-{p \over k}} (f_{k-1}(\Delta))^{{p \over k}},$$
which contradicts Theorem~\ref{equalvertices}.  \endproof

What this means in our context is that if we wish to choose the values of $a_j^p$ for some fixed $p$ and have already chosen the values of $a_j^k$ for all $k > p$, then we minimize the number of faces of smaller dimensions used at this step by setting $a_1^p = a_2^p = \dots = a_p^p$.  This is why the simple approach of making the constants $a_j^k$ not depend on $j$ is a reasonable first attempt.

However, we have seen that this approach sometimes fails.  The next conjecture explains why.

\begin{conjecture}  \label{unequalvertices}
Let $1 \leq p < k < r$ be integers.  Let $m$ and $n$ be real numbers such that $m > {r \choose k} n^{{k \over r}}$.  Define $w(i) = k$ if $i \leq k$ and $w(i) = r$ if $i > k$.  Then among all possible choices of positive real numbers $a_1^k, a_2^k, \dots a_k^k, a_1^r, a_2^r, \dots a_r^r$ that satisfy
\begin{enumerate}
\item $n = a_1^r a_2^r \dots a_r^r,$
\item $$m = \sum_{i_1 < i_2 < \dots < i_k} a_{i_1}^{w(i_k)} a_{i_2}^{w(i_k)} \dots a_{i_k}^{w(i_k)},$$
\item for every $i$, $j$, and $k$ with $i < j$, $a_i^k \geq a_j^k$, and
\item for every $i$, $j$, and $k$ with $i < k$, $a_j^i \geq a_j^k$,
\end{enumerate}
the unique choice of constants $a_j^k$ that minimizes
$$\sum_{i_1 < i_2 < \dots < i_p} a_{i_1}^{w(i_p)} a_{i_2}^{w(i_p)} \dots a_{i_p}^{w(i_p)}$$
is also the unique choice of constants $a_j^k$ such that
$$a_1^k = a_2^k = \dots = a_k^k = a_1^r = a_2^r = \dots = a_k^r > a_{k+1}^r = \dots = a_r^r.$$
\end{conjecture}

Intuitively, this means that if we have a given number of faces of dimensions $k-1$ and $r-1$, the way to minimize the number of faces of dimension $p-1$ is that, when allocating the vertices among various colors to use the faces of dimension $r-1$, we use more vertices among the first $k$ colors and fewer among the final $r-k$ colors, with a discrepancy exactly large enough to use all of the faces of dimension $k-1$ while creating the faces of dimension $r-1$.  This is what Example~\ref{unbalancedexample} did, using $r = 3$, $k = 2$, and $p = 1$.

In practice, this means that if choosing constants $a_j^k$ by the naive method and we run into a problem where $a_1^i < a_1^{i+1}$, the solution is to go back and make $a_j^k$ larger for some values of $j \leq i+1 < k$.  This necessarily means making $a_j^k$ smaller for some values of $i+1 < j \leq k$.

This is, of course, only a heuristic, and not a deterministic method.  A deterministic method for determining whether constants $a_j^k$ as in Conjecture~\ref{bigconj} exist, and giving them if they do, would surely be useful.

\section{$h$-vectors of vertex-decomposable flag complexes}

In this section, we explain how the constructions of this paper can be adapted to build a vertex-decomposable flag complex with a specified $h$-vector.  Furthermore, because every vertex-decomposable complex is constructible, shellable, and Cohen-Macaulay, we can replace ``vertex-decomposable" by any of these other terms and the results of this section will still hold.

The definitions of vertex-decomposable, constructible, shellable, and Cohen-Macaulay are rather technical, and not directly needed in this paper.  As such, we do not explicitly define them in this paper.  Loosely, they are conditions on simplicial complexes that ensure that the complex can be assembled nicely.

Given a simplicial complex $\Delta$ of dimension $d-1$, we can define its \textit{$h$-numbers} by
$$h_k(\Delta) = \sum_{i=0}^k (-1)^{k-i} {d-i \choose k-i} f_{i-1}(\Delta).$$
Conversely, we can compute
$$f_{k-1}(\Delta) = \sum_{i=0}^k {d-i \choose k-i} h_i(\Delta).$$
The \textit{$h$-vector} of $\Delta$ is the list of $h$-numbers, $h(\Delta) = (h_0, h_1, \dots, h_d)$.  The face vector and $h$-vector are merely different ways of presenting in the same information.  The face vector and $h$-vector are also related by the polynomial equation
$$\sum_{i=0}^d f_{i-1}(\Delta) t^{d-i} = \sum_{i=0}^d h_i(\Delta) (t+1)^{d-i}.$$

In some situations, it is more convenient to work with the $h$-vector of a simplicial complex rather than the face vector.

The \textit{chromatic number} of a simplicial complex is the minimum number of colors needed to color every vertex of the complex such that no two vertices in the same face are the same color.  This is equivalent to the chromatic number of the 1-skeleton of the complex, taken as a graph, in the usual graph-theoretic sense.  A complex is \textit{balanced} if its chromatic number is exactly one greater than its dimension.  A simplicial complex of dimension $d-1$ must have a face on $d$ vertices, all of which are distinct colors, so its chromatic number must be at least $d$.  For the complex to be balanced means that the chromatic number is not greater than $d$.

\begin{construction} \label{cmconst}
\textup{Let $\Delta$ be a colored complex.  Define a complex $\Delta^+$ as follows.  The vertices of $\Delta^+$ consist of all vertices of $\Delta$ together with one additional vertex of each color.  A set $F$ of vertices in $\Delta^+$ forms a face of $\Delta^+$ exactly if no two are the same color and the vertices of $F$ that are also in $\Delta$ form a face of $\Delta$.}  \endproof
\end{construction}

It is immediate from the construction that if $\Delta$ is a flag complex, then $\Delta^+$ is also a flag complex.  Furthermore, if $\Delta$ is balanced, then $\Delta^+$ is of the same dimension as $\Delta$.

\begin{theorem}[Cook-Nagel]  \label{dpshell}
Let $\Delta$ be a colored flag complex.  Then $\Delta^+$ is vertex-decomposable.
\end{theorem}

This was proven by Cook and Nagel in \cite[Theorem~3.3]{cook}.  Their terminology is very different from ours, so we should explain the connection.  The independence complex of a graph is formed by reversing the edges.  That is, two vertices form an edge in the independence complex exactly if they do not form an edge in the graph.  The independence complex is then the clique complex of the graph with its edges reversed.  If $\Delta$ is our colored flag complex, then the graph $G$ in \cite[Theorem~3.3]{cook} should be a graph with the edges of $\Delta$ reversed.  The vertex set $W_i$ in their clique-vertex partition is the set of vertices of color $i$ in our $\Delta$.  With this definition, their Ind $G^{\pi}$ is our $\Delta^+$.

\begin{proposition}  \label{dpfh}
Let $\Delta$ be a balanced, colored complex with color set $[n]$.  Then for every $0 \leq i \leq n$, $h_i(\Delta^+) = f_{i-1}(\Delta)$.
\end{proposition}

This result has been proven by many people and in a number of different ways.  For a recent example, see \cite[Proposition~3.8]{cook}.

Construction~\ref{mainconst} builds a flag complex of dimension $d-1$ that explicitly breaks most of the vertices into $d$ parts, with no two vertices in the same part adjacent.  We can give each part its own color.  The ``extra" vertices don't have a specified part, but their neighbors are all among at most $d-1$ of the other parts, so we can also give each of them whichever color we pointedly excluded from being used for any of their neighbors.  Therefore, any flag complex resulting from Construction~\ref{mainconst} is balanced by construction, even if we use the tweaks to make the number of vertices in distinct parts unequal.

Thus, if we want a vertex-decomposable flag complex with a specified $h$-vector, we can start by using Construction~\ref{mainconst} to construct a flag complex $\Delta$ with the specified face vector.  If we succeed, then we can add a vertex of each color to get a shellable flag complex $\Delta^+$ with the specified $h$-vector.

\end{document}